\documentclass[preprint,11pt,times]{elsarticle}

\usepackage{amsmath,amssymb,amsfonts,amsthm}
\usepackage{graphicx}
\usepackage{booktabs}
\usepackage{array}
\usepackage[colorlinks=true,linkcolor=blue,citecolor=blue,urlcolor=blue]{hyperref}
\usepackage{microtype}
\usepackage{bm}
\usepackage{mathtools}
\usepackage{xcolor}

\newtheorem{theorem}{Theorem}[section]
\newtheorem{lemma}[theorem]{Lemma}
\newtheorem{proposition}[theorem]{Proposition}
\newtheorem{corollary}[theorem]{Corollary}
\theoremstyle{definition}

\theoremstyle{remark}
\newtheorem{remark}{Remark}[section]

\newcommand{\e}{\mathrm{e}}
\newcommand{\vect}[1]{\bm{#1}}
\newcommand{\mat}[1]{\mathbf{#1}}

\newcommand{\Mtrans}{\mat{M}_{\mathrm{trans}}}

\DeclareMathOperator{\diag}{diag}
\DeclareMathOperator{\spr}{spr}

\graphicspath{{./}}
\journal{Journal of Computational Physics}

\AtBeginDocument{\hypersetup{%
  pdfauthor={Jin Zhao},
  pdftitle={Multiple-relaxation-time vectorial lattice Boltzmann schemes for incompressible flow: second-order consistency and parameter selection}%
}}

\begin{document}

\begin{frontmatter}

\title{Multiple-relaxation-time vectorial lattice Boltzmann schemes for incompressible flow: second-order consistency and parameter selection}

\author[inst0]{Yuyang Xu}
\ead{1200503040@cnu.edu.cn}

\affiliation[inst0]{organization={School of Mathematical Sciences, Capital Normal University},
            city={Beijing},
            postcode={100048},
            country={China}}

\author[inst1]{Yuxuan Zhang}

\author[inst1]{Jin Zhao\corref{cor1}}
\ead{zjin@cnu.edu.cn}
\cortext[cor1]{Corresponding author}

\affiliation[inst1]{organization={Academy for Multidisciplinary Studies, Capital Normal University},
            city={Beijing},
            postcode={100048},
            country={China}}

\begin{abstract}
We construct and analyze multiple-relaxation-time vectorial lattice Boltzmann
schemes for the incompressible Navier--Stokes equations in two and three space
dimensions.  The same moment-space relaxation matrix is applied to all physical
components, with independent rates for the flux, stress, and energy moments.
A Maxwell expansion organized by moment parity establishes formal second-order
consistency of the unforced interior schemes under diffusive scaling, assuming
a smooth near-equilibrium Maxwell branch.  The viscosity depends only on the
flux rate; for the two-dimensional scheme, we identify the first contributions
of the stress and energy rates to the second-order macroscopic truncation error.
We further construct a positive definite symmetrizer for the collision Jacobian
at rest and derive an exact weighted-norm identity for the nonlinear kinetic
residual under local collision.  At fixed viscosity, a finite Fourier search
selects the free rates by combining transverse-mode accuracy with
short-wavelength damping.  The resulting rates are compared with those
calibrated against flow errors on a coarse grid.  Numerical experiments exhibit
second-order convergence and illustrate how the free rates affect spectral
properties and finite-time flow errors.
\end{abstract}


\end{frontmatter}

\section{Introduction}
\label{sec:intro}

The lattice Boltzmann method (LBM) advances discrete populations by
alternating a local collision update with streaming on a prescribed velocity
stencil; see Refs.~\cite{Chen1998,Succi2001,Guo2013,Kruger2017}.  For a fixed
stencil, the relaxation spectrum of the collision operator affects the
dispersion, dissipation, and linear stability of the fully discrete scheme.
The Bhatnagar--Gross--Krook (BGK), or single-relaxation-time (SRT), model uses
one relaxation rate for all nonconserved modes~\cite{Bhatnagar1954}, whereas
two-relaxation-time (TRT) and MRT models assign rates to parity
classes or groups of moments~\cite{DHumieres1992,Lallemand2000,Ginzburg2008,Luo2011,Coreixas2019}.
Cascaded, cumulant, and entropic formulations provide alternative collision
constructions~\cite{Geier2006,Geier2015,Karlin1999,Karlin2014}.

Discrete kinetic approximations of incompressible flow based on vector BGK
models were developed by Carfora and Natalini~\cite{CarforaNatalini2008}.
Second-order BGK flux vector splitting schemes for the incompressible
Navier--Stokes equations were subsequently analyzed by Bouchut
\emph{et al.}~\cite{BouchutEtAl2018BGKFVS}.
The vectorial lattice Boltzmann method (VLBM) for incompressible flow uses vector-valued
populations whose components represent density and momentum
variables~\cite{ZhaoZhangYong2020Approximation,Zhao2021VectorBGK}.  Maxwell
iteration under diffusive scaling relates the discrete kinetic equations
to their macroscopic counterparts~\cite{zhao2017maxwell}, and weighted energy estimates have been
established for the associated discrete kinetic
approximation~\cite{ZhaoYong2020Weighted}.  An earlier D3N7 vectorial BGK
construction treats nonlinear pressure approximations and parameterized
second-order boundary conditions and derives a symmetrizer for its collision
operator~\cite{ZhaoZhang2020D3N7}.  A vectorial LBM based on finite differences has also
been analyzed for consistency and stability~\cite{ZhaoYong2024VFDLBM}.
A wall treatment for vector-valued distributions on a discrete velocity set was
introduced in Ref.~\cite{ZhaoZhangYong2020Boundary}; for the scheme considered
there, the midpoint rule is second-order accurate, while parameterized
extensions accommodate general boundary locations.

Relaxation has been extended in several ways within vectorial kinetic
formulations.  For nonlinear hyperbolic systems, Dubois introduced a
moment-space framework allowing a general MRT update; the shallow-water
examples reported there use equal relaxation rates~\cite{Dubois2014Vectorial}.
Dellar assigned separate relaxation times to the symmetric and antisymmetric
parts of the momentum flux to obtain a projection onto symmetric flux
tensors~\cite{Dellar2025VectorProjection}.  For incompressible flow,
Aregba-Driollet \emph{et al.} assigned distinct BGK-type relaxation rates to
the density and momentum components and analyzed second-order consistency
and linearized spectral properties~\cite{AregbaDriolletEtAl2026Vectorial}.
Within each physical component, all nonconserved velocity moments share one
rate.  For a common rate below unity, they also established entropy
dissipation.  For the vectorial TRT family based on the even and odd parts
associated with opposite links,
Ref.~\cite{Zhao2026TRTPower} establishes mesh-uniform power stability under a
complementary relation for homogeneous transport operators
that are unitary in the equilibrium metric and reversible under velocity
exchange.

Here we construct D2N5 and D3N7 vectorial MRT schemes in which a common
velocity-moment relaxation operator acts on every physical component.
Independent rates are assigned to the flux, stress, and energy moment rows.
We determine how this separation of moment relaxation preserves formal
second-order consistency under diffusive scaling, which rate determines the
viscosity, and where the remaining rates enter the conserved expansion.
The resulting freedom is then examined through parameter selection at fixed
viscosity.

For the equilibria considered here, making the relaxation matrix diagonal in
the chosen moment basis does not diagonalize the Jacobian of the nonlinear
collision map.  At the uniform rest state, the equilibrium flux moments vanish,
but their derivatives with respect to the conserved variables do not.
Moreover, a weighted contraction estimate for local collision does not by
itself control the full update unless streaming is nonexpansive in the same
norm.  We therefore analyze the local collision map separately from the
combined linearized collision and streaming operator for periodic problems.
For lattice Boltzmann schemes with scalar populations, asymptotic and Fourier analyses have clarified the
influence of nonhydrodynamic modes on hydrodynamic limits and short-wave
behavior~\cite{Dellar2002,Dellar2003IncompressibleLimits,Wissocq2019,Hosseini2020}.
Optimization over wavenumbers has also been used to determine free MRT
parameters~\cite{XuSagaut2011Optimal,XuMalaspinasSagaut2012Sensitivity}.
Yang \emph{et al.} established algebraic sufficient conditions for weighted
stability of linearized scalar-population MRT schemes using a diagonal
population-space metric~\cite{YangZhaoLin2024Automatic}.  That metric also
controls periodic streaming.  The generally non-diagonal symmetrizer
constructed here gives a local collision estimate; the combined update is
examined separately by Fourier analysis.

The analysis first identifies the rate relations under which the vectorial
MRT collision map reduces to TRT or SRT.  A Maxwell expansion organized by
moment parity establishes formal second-order consistency of the unforced
interior schemes under diffusive
scaling, assuming a smooth Maxwell branch with the stated uniform remainder
bounds.  Only the flux rate enters the viscosity, leaving the stress and
energy rates free at fixed viscosity.  For D2N5, the explicit
third-degree conserved coefficient identifies the first contribution of the
stress and energy rates to the $O(h^2)$ macroscopic truncation error.
As discussed after Eq.~\eqref{eq:free_rates_c3}, this is a partial contribution;
the fourth-degree conserved coefficient may contribute at the same order.
A lower triangular change of coordinates diagonalizes the collision Jacobian
at rest and yields a positive definite symmetrizer.
Weighted moment orthogonality gives an exact contraction identity for the
nonlinear kinetic residual under local collision, during which the conserved
variables and their equilibrium are unchanged.  These results describe the
interior consistency and local collision structure of the vectorial MRT
family; stability and convergence of the full initial-boundary-value problem
require additional analysis.

We use two distinct procedures for the remaining rates at fixed viscosity.
The first is a priori Fourier selection from the collision--streaming symbol
linearized about a uniform rest state.  It uses the same transverse
eigenvalue accuracy objective and short-wavelength damping criterion for MRT
and TRT, without fitting flow errors.  Finite wave-vector sampling and local
maximization provide numerical spectral checks; they do not establish a
uniform bound on the powers of the amplification matrix over all wave
vectors and time steps.  The second procedure is flow calibration.  For the
nonlinear D3 Taylor--Green (TGV) flow, we minimize coarse-grid velocity and
vorticity errors against an independent reference, then keep the rates fixed
in computations on finer grids and at a higher amplitude.  A separate channel-flow
calibration minimizes profile and flow-rate errors on one coarse Poiseuille
grid, with and without a spectral constraint.

Throughout the paper, \emph{Fourier-selected} MRT and TRT (F-MRT and F-TRT in
tables) refer only to the first procedure.  \emph{TGV-calibrated} and
\emph{channel-calibrated} rates identify the flow used by the second.
Complementary TRT (C-TRT) and SRT are reference choices fixed by their rate
relations, without either optimization procedure.

Periodic exact flows exhibit second-order numerical convergence.  The
nonlinear Taylor--Green tests assess the additional MRT freedom by retaining
coarse-grid calibrated rates on finer grids and at a higher amplitude.
These calibrated MRT rates reduce the errors relative to similarly
calibrated TRT, whereas the Fourier-selected MRT rates do not improve on
the TRT choices in this test.  Additional plane waves, a double shear layer,
and channel flows test how the rate choice affects different errors.  The
additional wave directions are absent from the transverse accuracy objective
but occur in the separate spectral test.  Complementary TRT gives the
smallest maximum error over these four waves at $n_g=32$, SRT gives the
smallest field errors in the double-shear tests, and channel flow-rate and
profile errors favor different choices.  Together, these comparisons
distinguish the effects of the free moment rates on the Fourier spectrum
and on finite-time flow accuracy.

Sections~\ref{sec:vlbm}--\ref{sec:d3n7} develop the D2N5 and D3N7 formulations
together with their consistency and collision analyses.
Sections~\ref{sec:algorithm} and~\ref{sec:numerical} present the Fourier
parameter selection and numerical experiments, respectively.

\section{Vectorial lattice Boltzmann discretization}
\label{sec:vlbm}

We first specify the D2N5 discretization.  The D3N7 analogue is given in
Section~\ref{sec:d3n7}.

\subsection{Discrete velocity stencil}

For $d=2$ and $N=5$, the discrete velocities are
\begin{equation}
    e_1 = \begin{pmatrix} 1 \\ 0 \end{pmatrix}, \quad
    e_2 = \begin{pmatrix} 0 \\ 1 \end{pmatrix}, \quad
    e_3 = \begin{pmatrix} -1 \\ 0 \end{pmatrix}, \quad
    e_4 = \begin{pmatrix} 0 \\ -1 \end{pmatrix}, \quad
    e_5 = \begin{pmatrix} 0 \\ 0 \end{pmatrix}.
    \label{eq:velocities}
\end{equation}
The opposite-direction involution $i\mapsto\bar i$ is defined by
$e_{\bar i}=-e_i$, with
$(\bar1,\bar2,\bar3,\bar4,\bar5)=(3,4,1,2,5)$.  More generally, the axis-aligned
$D_dN_{2d+1}$ velocity set contains the $2d$ signed coordinate vectors and the rest
velocity.

\subsection{Vector-valued populations}

For each velocity $e_i$, let the corresponding population be
\begin{equation}
    f_i = \begin{pmatrix} f_{i,0} \\ f_{i,1} \\ \vdots \\ f_{i,d} \end{pmatrix}
    \in \mathbb{R}^{d+1}, \qquad i=1,\ldots,N.
    \label{eq:vector_distribution}
\end{equation}
The discrete state is the stacked vector
\begin{equation}
    \vect{f} = \begin{pmatrix} f_1 \\ f_2 \\ \vdots \\ f_N \end{pmatrix}
    \in \mathbb{R}^{N(d+1)}.
\end{equation}
Thus $\vect f\in\mathbb R^{15}$ for D2N5.

\subsection{Conserved variables}

The conserved moment is
\begin{equation}
    w = \sum_{i=1}^{N} f_i = \begin{pmatrix} \rho \\ h\rho u_1 \\ \vdots \\ h\rho u_d \end{pmatrix} \in \mathbb{R}^{d+1},
    \label{eq:conserved}
\end{equation}
where $\rho$ denotes density, $\vect u=(u_1,\ldots,u_d)^T$ denotes velocity,
and $h=\Delta x$.  The factor $h$ in the momentum variables is consistent with
the diffusive scaling.

\subsection{Equilibrium map}

Let $a\in(0,1/4)$ and $\Delta t=\alpha h^2$.  The D2N5 equilibrium
populations are defined by
\begin{align}
    M_1 &= aw + \frac{\alpha}{2}A_1, \quad M_3 = aw - \frac{\alpha}{2}A_1, \label{eq:eq_dist_1}\\
    M_2 &= aw + \frac{\alpha}{2}A_2, \quad M_4 = aw - \frac{\alpha}{2}A_2, \label{eq:eq_dist_2}\\
    M_5 &= (1-4a)w, \label{eq:eq_dist_5}
\end{align}
with flux vectors
\begin{equation}
    A_1 = \begin{pmatrix} h\rho u_1 \\ h^2\rho u_1^2 + h^2 P \\ h^2\rho u_1 u_2 \end{pmatrix}, \quad
    A_2 = \begin{pmatrix} h\rho u_2 \\ h^2\rho u_1 u_2 \\ h^2\rho u_2^2 + h^2 P \end{pmatrix},
    \label{eq:flux_vectors}
\end{equation}
where $P=(\rho-\bar\rho)/h^2$ is the pressure relative to the reference
density $\bar\rho$.

The equilibrium has the same conserved moment as $\vect f$, since
\begin{equation}
    \sum_{i=1}^{N} M_i = 4aw + (1-4a)w = w.
\end{equation}

\section{MRT collision operator}
\label{sec:mrt}

\subsection{Moment-space notation}

Let $\Mtrans\in\mathbb R^{N\times N}$ act on the velocity index and let
$\mat S=\diag(s_0,s_1,\ldots,s_{N-1})$.  The corresponding operators on the
stacked vector population are
\begin{equation}
    \tilde{\mat{M}} = \Mtrans \otimes I_{d+1}, \quad \tilde{\mat{S}} = \mat{S} \otimes I_{d+1}.
\end{equation}
For D2N5, both operators have dimension $15\times15$.

\subsection{Collision--streaming discretization}

One time step is defined by
\begin{equation}
    f_i(\vect{x} + e_i h, t + \alpha h^2) = f_i(\vect{x}, t) + Q_i(\vect{f}) + \alpha h^3 F_i,
    \label{eq:evolution}
\end{equation}
where $F_i=0$ for $i\ne N$ and $F_N=F$ is the prescribed body force vector.
The local collision operator is
\begin{equation}
    \vect{Q}(\vect{f}) = -\tilde{\mat{M}}^{-1}\tilde{\mat{S}}\tilde{\mat{M}}(\vect{f} - \vect{f}^{\rm eq}),
    \label{eq:collision}
\end{equation}
with $\vect f^{\rm eq}=(M_1,\ldots,M_N)^T$.
Because the first moment is conserved and
$\sum_iM_i=w=\sum_if_i$, Eq.~\eqref{eq:collision} gives
$\sum_iQ_i=0$; hence collision preserves $w$.

For $\vect m=\tilde{\mat M}\vect f$ and
$\vect m^{\rm eq}=\tilde{\mat M}\vect f^{\rm eq}$, the collision step reads
\begin{equation}
    m_k' = m_k - s_k(m_k - m_k^{\rm eq}), \quad k = 0, 1, \ldots, N-1,
\end{equation}
where $m_k\in\mathbb R^{d+1}$.  Thus the prescribed relaxation matrix is
diagonal in the raw moment coordinates; the full collision Jacobian need not be
diagonal because $\vect m^{\rm eq}$ depends on the conserved moment.

\subsection{Weighted orthogonality of the moment basis}

Let
\begin{equation}
    \mat W=\diag(a,a,a,a,1-4a)
    \label{eq:d2_direction_weights}
\end{equation}
be the positive diagonal matrix of direction weights, and define
$\langle r,q\rangle_{\mat W}=r\mat Wq^T$ for row vectors $r$ and $q$.  We use a
$\mat W$-orthogonal moment basis whose first row represents the conserved
moment and whose stress and energy rows have zero equilibrium values.  With
$\beta=4a/(1-4a)$, such a basis is
\begin{equation}
    \Mtrans = \begin{pmatrix}
    1 & 1 & 1 & 1 & 1 \\
    1 & 0 & -1 & 0 & 0 \\
    0 & 1 & 0 & -1 & 0 \\
    1 & -1 & 1 & -1 & 0 \\
    1 & 1 & 1 & 1 & -\beta
    \end{pmatrix}.
    \label{eq:M_trans}
\end{equation}

The corresponding raw moments and their equilibrium values are
\begin{equation}
\begin{aligned}
    m_0 &= \sum_{i=1}^{5}f_i,
        & m_0^{\rm eq} &= w,\\
    m_1 &= f_1-f_3,
        & m_1^{\rm eq} &= \alpha A_1,\\
    m_2 &= f_2-f_4,
        & m_2^{\rm eq} &= \alpha A_2,\\
    m_3 &= f_1-f_2+f_3-f_4,
        & m_3^{\rm eq} &= 0,\\
    m_4 &= f_1+f_2+f_3+f_4-\beta f_5,
        & m_4^{\rm eq} &= 0.
\end{aligned}
\label{eq:d2_raw_moments}
\end{equation}
Thus $m_0$ is conserved, $m_1$ and $m_2$ are the flux moments, and $m_3$ and
$m_4$ are the even stress and energy moments, respectively.  The equilibrium
values in~\eqref{eq:d2_raw_moments} follow from
Eqs.~\eqref{eq:eq_dist_1}--\eqref{eq:eq_dist_5}.  Here ``stress'' and
``energy'' label rows of the discrete-velocity moment basis; they do not
introduce a separate constitutive stress or physical energy equation.

The inverse matrix appearing in~\eqref{eq:collision} is
\begin{equation}
    \Mtrans^{-1} = \begin{pmatrix}
    a & \frac{1}{2} & 0 & \frac{1}{4} & \frac{1-4a}{4} \\[4pt]
    a & 0 & \frac{1}{2} & -\frac{1}{4} & \frac{1-4a}{4} \\[4pt]
    a & -\frac{1}{2} & 0 & \frac{1}{4} & \frac{1-4a}{4} \\[4pt]
    a & 0 & -\frac{1}{2} & -\frac{1}{4} & \frac{1-4a}{4} \\[4pt]
    1-4a & 0 & 0 & 0 & -(1-4a)
    \end{pmatrix}.
    \label{eq:M_inv}
\end{equation}

\subsection{TRT restriction of the MRT collision map}

A TRT collision operator assigns one relaxation rate to the nonconserved even
subspace and one to the odd subspace.  In terms of parameters $(\phi,\psi)$,
the update is
\begin{equation}
    f_i' = f_i + \phi(M_i - f_i) + \psi(M_{\bar{i}} - f_{\bar{i}}).
\end{equation}

Let $\mat P$ denote the opposite-direction permutation,
$P_{ij}=\delta_{j,\bar i}$.  Then
\begin{equation}\label{eq:trt_matrix}
    \vect{f}' = \vect{f}
    + \bigl[(\phi \mat{I} + \psi \mat{P})\otimes I_3\bigr]
      (\vect{f}^{\rm eq} - \vect{f}).
\end{equation}

Define the even and odd parts by
\begin{equation}
    f_i^\pm = \frac{f_i \pm f_{\bar{i}}}{2}, \qquad
    M_i^\pm = \frac{M_i \pm M_{\bar{i}}}{2}.
\end{equation}
Equivalently,
\begin{equation}
    f_i' = f_i - s^+(f_i^+ - M_i^+) - s^-(f_i^- - M_i^-),
\end{equation}
where $s^+=\phi+\psi$ and $s^-=\phi-\psi$ are the relaxation rates of the even
and odd subspaces, respectively.

The flux moments \(m_1=2f_1^-\) and \(m_2=2f_2^-\) are odd under
opposite-velocity exchange, whereas
\[
 m_3=2f_1^+-2f_2^+,
 \qquad m_4=2f_1^++2f_2^+-\beta f_5
\]
and the conserved moment \(m_0\) are even.

\begin{lemma}[Parity decomposition in moment coordinates]\label{lem:P_diag}
The opposite-direction permutation is diagonal in the chosen moment basis:
\begin{equation}
    \Mtrans \mat{P} \Mtrans^{-1} = \mat{D}_{P}, \qquad
    \mat{D}_{P} := \diag(1, -1, -1, 1, 1).
\end{equation}
\end{lemma}

\begin{proof}
The rows of \(\Mtrans\) are left eigenvectors of \(\mat P\).  The parity
identified above gives the eigenvalues \((1,-1,-1,1,1)\), respectively.
\end{proof}

\begin{theorem}[TRT representation of the MRT collision map]\label{thm:mrt_trt}
Let $\Mtrans$ be an invertible moment matrix whose first row is
$\vect 1_5^T$, representing the conserved moment.  Suppose its rows are
left eigenvectors of the opposite-velocity permutation $\mat P$, ordered so that
\[
\Mtrans\mat P\Mtrans^{-1}
=\diag(1,-1,-1,1,1).
\]
For the equilibrium map defined above, the MRT and TRT collision maps coincide
for every state with positive density if and only if
\begin{align}
    s_1 = s_2 &= s^- = \phi - \psi, \label{eq:s_odd}\\
    s_3 = s_4 &= s^+ = \phi + \psi. \label{eq:s_even}
\end{align}
Equivalently,
\begin{equation}\label{eq:S_star}
 \mat S=\diag(s_0,s^-,s^-,s^+,s^+),
 \qquad s^-=\phi-\psi,\quad s^+=\phi+\psi,
\end{equation}
where $s_0$ is arbitrary.  Here $s^-$ acts on the two flux moments, whereas
$s^+$ acts on the stress and energy moments.
Conversely, imposing these pairwise equalities on the MRT rates determines
\begin{equation}\label{eq:phi_psi_from_mrt}
    \phi = \frac{s_1 + s_3}{2}, \quad \psi = \frac{s_3 - s_1}{2}.
\end{equation}
The underlying matrices $\Mtrans^{-1}\mat S\Mtrans$ and
$\phi\mat I+\psi\mat P$ are equal if and only if, in addition, $s_0=s^+$; in
that case,
\begin{equation}\label{eq:trt_mrt_matrix}
 \Mtrans^{-1}\mat S\Mtrans=\phi\mat I+\psi\mat P.
\end{equation}
\end{theorem}

\begin{proof}
By Lemma~\ref{lem:P_diag},
\begin{align*}
    \Mtrans (\phi \mat{I} + \psi \mat{P}) \Mtrans^{-1} &= \phi \mat{I} + \psi \mat{D}_{P} \\
    &= \diag(\phi + \psi, \phi - \psi, \phi - \psi, \phi + \psi, \phi + \psi) \\
    &= \diag(s^+, s^-, s^-, s^+, s^+).
\end{align*}
For the moment residual
$\vect y=\tilde{\mat M}(\vect f-\vect f^{\rm eq})$, conservation gives
$y_0=0$.
The nonconserved components of $\vect y$ can otherwise be chosen independently.
The two collision maps therefore coincide for every state with positive density if and
only if their diagonal actions agree in moment rows $1,\ldots,4$, which gives
Eqs.~\eqref{eq:s_odd}--\eqref{eq:s_even}.  The possible difference in the
conserved row annihilates the residual.  Equality of the underlying matrices also
requires agreement in row zero and hence $s_0=s^+$.
\end{proof}

The SRT collision map is the special case \(s^-=s^+\), equivalently
\(s_1=s_2=s_3=s_4\).  The parameter \(s_0\) in the conserved row remains immaterial
because the conserved collision residual vanishes identically.  Write
\(s_\nu=s_1=s_2\), \(s_\sigma=s_3\), and \(s_\varepsilon=s_4\), and let
\(\mathcal I=(0,2)\).  The corresponding sets of nonconserved rates are
\begin{equation}
 \begin{aligned}
 \mathcal P_{\rm SRT}&=\{(s,s,s):s\in\mathcal I\},\\
 \mathcal P_{\rm TRT}&=\{(s^-,s^+,s^+):(s^-,s^+)\in\mathcal I^2\},\\
 \mathcal P_{\rm MRT}&=\mathcal I^3.
 \end{aligned}
 \qquad
 \mathcal P_{\rm SRT}\subsetneq\mathcal P_{\rm TRT}
 \subsetneq\mathcal P_{\rm MRT}.
 \label{eq:srt_trt_mrt_nesting}
\end{equation}
For the present equilibrium and moment basis, these inclusions also hold
for the associated families of local collision maps.  They do not compare
the consistency, stability, or accuracy of the resulting schemes.  With
$s_\nu$ fixed by the viscosity, SRT gives a single point, TRT the line
$s_\sigma=s_\varepsilon$, and MRT the two-dimensional set
$\{s_\nu\}\times\mathcal I^2$.

\begin{remark}[Sign of the energy row]
Parity alone does not determine the sign of the rest-velocity coefficient in
the energy row: either sign leaves the row even under $\mat P$, provided the
moment matrix remains invertible.  The negative choice is required here to
obtain both $m_4^{\rm eq}=0$ and weighted orthogonality to the conserved row.
If $r_k$ denotes row $k$ of $\Mtrans$, then
\begin{equation*}
    r_0 \mat W r_4^T
    =4a-(1-4a)\frac{4a}{1-4a}=0.
\end{equation*}
With the positive sign, this inner product equals $8a\ne0$ and the energy
equilibrium is nonzero, although the MRT--TRT identity based on parity remains
valid for the grouped rates in~\eqref{eq:S_star}.
\end{remark}

We use the conventional value $s_0=0$ in the remainder of the paper.

\section{Formal second-order consistency and viscosity under diffusive scaling}
\label{sec:consistency}

Let \(h=\Delta x\) and \(\Delta t=\alpha h^2\), where \(\alpha>0\) is fixed, and
consider an \(h\)-independent time interval.  We use the weakly compressible
scaling
\begin{equation}
 w=(\rho,h\rho u_1,h\rho u_2)^T,
 \qquad \rho=\bar\rho+h^2P,
 \label{eq:cons_scaling}
\end{equation}
with \(\bar\rho>0\).  The fields \(P\) and \(\vect u\), and all derivatives
needed below, are assumed uniformly bounded.  Assume that the streaming Taylor
expansion through parabolic degree four has a componentwise \(O(h^5)\)
remainder, uniformly on the time interval.  Set
\(s_1=s_2=s_\nu\), \(s_3=s_\sigma\), and \(s_4=s_\varepsilon\), where the
rates are fixed as \(h\to0\) and
\[
 0<s_{\min}\le s_\nu<2,
 \qquad s_{\min}\le s_\sigma,s_\varepsilon\le2.
\]
For the Maxwell iteration, assume that a smooth near-equilibrium Maxwell
branch exists through parabolic degree four.  Introduce an auxiliary
parameter \(\lambda\), assigning degree one to a space derivative and degree
two to a time derivative, and write
\begin{equation}
 \vect m(\lambda)=\vect m^{\rm eq}
   +\sum_{r=1}^{4}\lambda^rh^r\vect m^{(r)}
   +\vect{\mathcal R}_m^{(5)}(\lambda;h),
 \qquad m_0^{(r)}=0,\quad (\mathcal R_m^{(5)})_0=0,
 \label{eq:cons_maxwell_ansatz}
\end{equation}
where the coefficients and their required derivatives are uniformly bounded,
and \(\vect{\mathcal R}_m^{(5)}\), together with the finitely many derivatives
used below, is componentwise \(O(h^5)\), uniformly for
\(0\le\lambda\le1\).  These are asymptotic regularity assumptions on the
Maxwell expansion; the conserved equations are not assumed.  The
coefficients are determined from the nonconserved projection below.
We set \(\lambda=1\) before imposing the conserved balance, rather than
matching that balance coefficientwise in \(\lambda\), because the physical
components in Eq.~\eqref{eq:cons_scaling} have different \(h\)-dependent
amplitudes.

\begin{proposition}[Formal second-order macroscopic consistency]
\label{prop:formal_consistency_scope}
For the smooth, unforced interior D2N5 scheme under the preceding assumptions,
the equivalent equations obtained by Maxwell iteration formally satisfy,
componentwise and with remainders uniform on the stated time interval,
\begin{align}
 \partial_t\rho+\nabla\!\cdot(\rho\vect u)&=O(h^2),
 \label{eq:cons_continuity_formal}\\
 \partial_t(\rho\vect u)
 +\nabla\!\cdot(\rho\vect u\otimes\vect u)
 +\nabla P-\nu\Delta(\rho\vect u)&=O(h^2),
 \label{eq:cons_momentum_formal}
\end{align}
where
\begin{equation}
 \nu=\frac{2a}{\alpha}\left(\frac1{s_\nu}-\frac12\right).
 \label{eq:viscosity_formula_precise}
\end{equation}
The conclusion holds for arbitrary fixed \(s_\sigma\) and \(s_\varepsilon\)
in the stated interval, including the TRT and SRT restrictions.  This formal
result for the equivalent equations does not establish stability or
convergence for a fully discrete initial-boundary-value problem.
\end{proposition}

\begin{proof}
From Eqs.~\eqref{eq:eq_dist_1}--\eqref{eq:eq_dist_5},
\begin{equation}
 \vect m^{\rm eq}=(w,\alpha A_1,\alpha A_2,0,0)^T.
 \label{eq:cons_meq}
\end{equation}
Define
\begin{equation}
 \widehat{\mat D}
 =(\Mtrans\otimes I_3)
 \left[\diag(\partial_x,\partial_y,-\partial_x,-\partial_y,0)
       \otimes I_3\right]
 (\Mtrans^{-1}\otimes I_3).
 \label{eq:cons_Dhat_definition}
\end{equation}
Its coefficients are constant, so \(\widehat{\mat D}\) commutes with
\(\partial_t\).  Direct multiplication gives
\begin{equation}
 (\widehat{\mat D}\vect m)_0=\partial_xm_1+\partial_ym_2,
 \quad
 (\widehat{\mat D}\vect m^{\rm eq})_{1,2}
   =2a(\partial_xw,\partial_yw),
 \quad
 (\widehat{\mat D}^{\,2}\vect m^{\rm eq})_0=2a\Delta w.
 \label{eq:cons_Dhat_identities}
\end{equation}
Set
\begin{equation}
 \widehat{\mat S}
 =\diag(0,s_\nu,s_\nu,s_\sigma,s_\varepsilon)\otimes I_3,
 \qquad
 \widehat{\mat S}^{\dagger}
 =\diag(0,s_\nu^{-1},s_\nu^{-1},s_\sigma^{-1},s_\varepsilon^{-1})
   \otimes I_3.
\end{equation}
Let \(\vect e_0\) be the first coordinate vector in \(\mathbb R^5\), and set
\[
 \mat P_c=(\vect e_0\vect e_0^T)\otimes I_3,
 \qquad \mat P_k=I_{15}-\mat P_c.
\]

With
\(\mathcal L_\lambda=\lambda h\widehat{\mat D}
+\lambda^2\alpha h^2\partial_t\), the Maxwell equation for the
nonconserved moments is
\begin{equation}
 \mat P_k\left\{[\exp(\mathcal L_\lambda)-I]\vect m(\lambda)
 +\widehat{\mat S}\bigl(\vect m(\lambda)-\vect m^{\rm eq}\bigr)\right\}=0.
 \label{eq:cons_taylor_operator}
\end{equation}
Put \(\vect m^{(0)}=\vect m^{\rm eq}\).  Write the coefficient of
\(\lambda^rh^r\) in the streaming term as
\begin{equation}
 \vect g^{(r)}=
 \sum_{\substack{p,q,\ell\ge0\\p+2q+\ell=r\\p+q\ge1}}
 \frac{\alpha^q}{p!q!}\,
 \widehat{\mat D}^{\,p}\partial_t^q\vect m^{(\ell)}.
 \label{eq:cons_maxwell_recursion}
\end{equation}
Because \(p+q\ge1\), every index in the sum satisfies \(\ell\le r-1\).
Hence \(\vect g^{(r)}\) depends only on
\(\vect m^{(0)},\ldots,\vect m^{(r-1)}\).  The nonconserved equations thus
determine \(\vect m^{(r)}\) successively, without using the conserved balance.
Matching the nonconserved block coefficientwise, for \(1\le r\le4\), gives
\begin{equation}
 \vect m^{(r)}=-\widehat{\mat S}^{\dagger}\vect g^{(r)},
 \qquad
 \vect c^{(r)}:=(\vect e_0^T\otimes I_3)\vect g^{(r)}.
 \label{eq:cons_recursion_compact}
\end{equation}

To control the discarded terms, let
\[
 \mat J=\diag(1,-1,-1,1,1)\otimes I_3,
 \qquad
 \mat\Pi_\tau=\frac12(I_{15}+\tau\mat J),
 \quad \tau\in\{+1,-1\}.
\]
Opposite-direction symmetry gives
\begin{equation}
 \mat J\widehat{\mat D}=-\widehat{\mat D}\mat J,
 \qquad [\mat J,\widehat{\mat S}^{\dagger}]=0.
 \label{eq:cons_Dhat_parity}
\end{equation}
The second relation does not require \(s_\sigma=s_\varepsilon\).

For \(\tau\in\{+1,-1\}\), let
\((\vect m_\tau^{(r)},\vect g_\tau^{(r)})\) denote the contribution to the
linear recursion
\eqref{eq:cons_maxwell_recursion}--\eqref{eq:cons_recursion_compact}
with the zeroth-order term
\(\vect m_\tau^{(0)}=\mat\Pi_\tau\vect m^{\rm eq}\).
By linearity,
\(\vect m^{(r)}=\vect m_+^{(r)}+\vect m_-^{(r)}\), and similarly for
\(\vect g^{(r)}\).  We claim that
\begin{equation}
 \mat J\vect m_\tau^{(r)}
   =\tau(-1)^r\vect m_\tau^{(r)}\quad(r\geq0),
 \qquad
 \mat J\vect g_\tau^{(r)}
   =\tau(-1)^r\vect g_\tau^{(r)}\quad(r\geq1).
 \label{eq:cons_projected_parity_induction}
\end{equation}
The assertion for \(\vect m_\tau^{(0)}\) follows from its definition.  If it
holds through order \(r-1\),
then a summand
\(\widehat{\mat D}^{\,p}\partial_t^q\vect m_\tau^{(\ell)}\)
in \(\vect g_\tau^{(r)}\) has parity
\[
 \tau(-1)^{p+\ell}=\tau(-1)^r,
\]
because \(p+2q+\ell=r\).  Equation
\eqref{eq:cons_Dhat_parity} and the recursion for
\(\vect m_\tau^{(r)}\) complete the induction.

Finally,
\[
 (\vect e_0^T\otimes I_3)\mat\Pi_-=0.
\]
Consequently, the conserved coefficient \(\vect c^{(r)}\) receives a
contribution from \(\mat\Pi_-\vect m^{\rm eq}\) only when \(r\) is odd, and
from \(\mat\Pi_+\vect m^{\rm eq}\) only when \(r\) is even.  Thus
\(\vect c^{(3)}\) is generated by the odd equilibrium part, whereas
\(\vect c^{(4)}\) is generated by the even equilibrium part.

The first two coefficients determine the macroscopic equation.
Equations~\eqref{eq:cons_Dhat_identities} and~\eqref{eq:cons_recursion_compact}
give
\begin{equation}
\begin{aligned}
 m_1^{(1)}&=-\frac{2a}{s_\nu}\partial_xw,\qquad
 m_2^{(1)}=-\frac{2a}{s_\nu}\partial_yw,\\
 \vect c^{(1)}&=\alpha(\partial_xA_1+\partial_yA_2),\\
 \vect c^{(2)}
 &=\alpha\partial_tw+(\widehat{\mat D}\vect m^{(1)})_0
   +\tfrac12(\widehat{\mat D}^{\,2}\vect m^{\rm eq})_0
 =\alpha(\partial_tw-\nu\Delta w),
\end{aligned}
\label{eq:cons_leading_coefficients}
\end{equation}
where \(\nu\) is given by Eq.~\eqref{eq:viscosity_formula_precise}.  After
setting \(\lambda=1\), the identity \(\mat P_c\widehat{\mat S}=0\) gives the
conserved projection of the formally reconstructed update:
\begin{equation}
 0=(\vect e_0^T\otimes I_3)
   [\exp(\mathcal L_1)-I]\vect m(1)
  =\sum_{r=1}^{4}h^r\vect c^{(r)}+\vect{\mathcal R}^{(5)}.
 \label{eq:cons_physical_projection}
\end{equation}
Here \(\vect{\mathcal R}^{(5)}\) collects all contributions of parabolic
degree at least five.  The boundedness and Taylor remainder assumptions make it
componentwise \(O(h^5)\).  Consequently,
\begin{equation}
 0=\alpha h(\partial_xA_1+\partial_yA_2)
   +\alpha h^2(\partial_tw-\nu\Delta w)+\vect{\mathcal R}_c,
 \label{eq:cons_conserved_expansion}
\end{equation}
where
\(\vect{\mathcal R}_c=h^3\vect c^{(3)}+h^4\vect c^{(4)}
+\vect{\mathcal R}^{(5)}\).

It remains to estimate \(\vect{\mathcal R}_c\).  Under
Eq.~\eqref{eq:cons_scaling}, the orders of the density and momentum components
in the odd and nonconstant even equilibrium parts are, respectively,
\begin{equation}
 (O(h),O(h^2),O(h^2)),
 \qquad (O(h^2),O(h),O(h)).
 \label{eq:cons_equilibrium_component_orders}
\end{equation}
Derivatives preserve these orders.  The parity conclusion above then yields
\[
 \vect c^{(3)}=(O(h),O(h^2),O(h^2))^T,
 \qquad
 \vect c^{(4)}=(O(h^2),O(h),O(h))^T.
\]
The constant part of the reference density contributes nothing because every term in
Eq.~\eqref{eq:cons_maxwell_recursion} contains a derivative.  Multiplication by
the parabolic factors \(h^3\) and \(h^4\), together with the fifth-degree
remainder, gives
\begin{equation}
 \vect{\mathcal R}_c=(O(h^4),O(h^5),O(h^5))^T.
 \label{eq:cons_conserved_remainder}
\end{equation}

Substituting Eqs.~\eqref{eq:cons_scaling} and~\eqref{eq:flux_vectors} directly
into Eq.~\eqref{eq:cons_conserved_expansion}, and using
Eq.~\eqref{eq:cons_conserved_remainder}, gives
\begin{align*}
 \alpha h^2[\partial_t\rho+\nabla\!\cdot(\rho\vect u)-\nu\Delta\rho]
   &=O(h^4),\\
 \alpha h^3[\partial_t(\rho u_j)+\partial_k(\rho u_k u_j)
       +\partial_jP-\nu\Delta(\rho u_j)]&=O(h^5),
 \qquad j=1,2,
\end{align*}
where the repeated index \(k\) is summed.  The scaling in
Eq.~\eqref{eq:cons_scaling} gives the prefactors \(\alpha h^2\) and
\(\alpha h^3\); division by these factors yields \(O(h^2)\) residuals in
both equations.  Since \(\Delta\rho=h^2\Delta P\), the density diffusion
term is absorbed into the continuity remainder.  This proves
Eqs.~\eqref{eq:cons_continuity_formal}--\eqref{eq:cons_momentum_formal}.  The
continuity equation and \(\rho=\bar\rho+h^2P\) also give
\(\nabla\!\cdot\vect u=O(h^2)\).  Finally,
Eq.~\eqref{eq:cons_leading_coefficients} contains neither \(s_\sigma\) nor
\(s_\varepsilon\), so these rates do not alter the displayed leading-order
macroscopic coefficients.  They may nevertheless enter the \(O(h^2)\)
macroscopic truncation term and hence its error constants.
\end{proof}

\begin{remark}[First occurrence of the free even relaxation rates]
\label{rem:free_rates_first_occurrence}
The stress and energy rates do not enter
\(\vect c^{(1)}\) or \(\vect c^{(2)}\), but their first-order kinetic
corrections are
\[
 m_3^{(1)}
   =-\frac{\alpha}{s_\sigma}
      (\partial_xA_1-\partial_yA_2),
 \qquad
 m_4^{(1)}
   =-\frac{\alpha}{s_\varepsilon}
      (\partial_xA_1+\partial_yA_2).
\]
Writing
\([\vect c^{(3)}]_{\sigma,\varepsilon}\) for the sum of the terms in
\(\vect c^{(3)}\) that contain \(s_\sigma^{-1}\) or
\(s_\varepsilon^{-1}\), direct substitution in the \(r=2\) and \(r=3\)
recursions gives
\begin{equation}
\begin{aligned}
 [\vect c^{(3)}]_{\sigma,\varepsilon}
  ={}&\frac{\alpha}{2}
       \left(\frac1{s_\nu}-\frac12\right)
       \bigg[
       \frac{1}{s_\sigma}
       (\partial_x^2-\partial_y^2)
       (\partial_xA_1-\partial_yA_2)\\
     &\hspace{42mm}
       +\frac{1-4a}{s_\varepsilon}
       \Delta(\partial_xA_1+\partial_yA_2)
       \bigg].
\end{aligned}
\label{eq:free_rates_c3}
\end{equation}
Thus the free even rates first enter the conserved expansion through the
successive couplings
\[
 \mat\Pi_-\vect m^{\rm eq}
 \longrightarrow \mat\Pi_+\vect m^{(1)}
 \longrightarrow \mat\Pi_-\vect m^{(2)}
 \longrightarrow \mat\Pi_+\vect g^{(3)}.
\]
Under the weakly compressible scaling, the density and momentum components
of Eq.~\eqref{eq:free_rates_c3} are \(O(h)\) and \(O(h^2)\), respectively.
After multiplication by \(h^3\) and division by the corresponding physical
prefactors, they contribute only to the \(O(h^2)\) macroscopic truncation
term.  Thus \(s_\sigma\) and \(s_\varepsilon\) can affect the error
constants, but Eq.~\eqref{eq:free_rates_c3} does not imply that they reduce
the error.  The displayed differential combinations may vanish for
particular flows.  Moreover, $\vect c^{(4)}$ may also contribute at the
same macroscopic order, so Eq.~\eqref{eq:free_rates_c3} gives only part of
the $O(h^2)$ truncation term.
\end{remark}

\section{Collision structure in modified moment coordinates}
\label{sec:stability}

We linearize the collision operator about a uniform rest state.
Although the equilibrium fluxes vanish at this state, their derivatives
couple density and momentum.

\subsection{Equilibrium Jacobian at rest}

Let
\[
    \mat L=\vect 1_5^T\otimes I_3,
    \qquad w=\mat L\vect f=(\rho,q_1,q_2)^T,
    \qquad q_k=h\rho u_k,
\]
and write the equilibrium map as
\(\vect f^{\rm eq}=\mathcal M(w)\).  In the conserved variables,
\begin{equation}
\begin{aligned}
 A_1(w)&=\left(q_1,\;q_1^2/\rho+\rho-\bar\rho,\;
                    q_1q_2/\rho\right)^T,\\
 A_2(w)&=\left(q_2,\;q_1q_2/\rho,\;
                    q_2^2/\rho+\rho-\bar\rho\right)^T.
\end{aligned}
\label{eq:flux_conserved_variables}
\end{equation}
At the uniform rest state \(w_*=(\bar\rho,0,0)^T\), with
\(\vect f_*=\mathcal M(w_*)\), define
\begin{equation}
 B_1=D_wA_1(w_*)=
 \begin{pmatrix}0&1&0\\1&0&0\\0&0&0\end{pmatrix},
 \qquad
 B_2=D_wA_2(w_*)=
 \begin{pmatrix}0&0&1\\0&0&0\\1&0&0\end{pmatrix}.
\label{eq:flux_jacobian_blocks}
\end{equation}

\begin{lemma}[Full D2N5 equilibrium Jacobian at rest]
\label{lem:full_equilibrium_jacobian}
Set
\[
\begin{array}{lll}
 K_1=aI_3+\frac{\alpha}{2}B_1,&
 K_2=aI_3+\frac{\alpha}{2}B_2,&
 K_3=aI_3-\frac{\alpha}{2}B_1,\\[2pt]
 K_4=aI_3-\frac{\alpha}{2}B_2,&
 K_5=(1-4a)I_3,&
 \mat K=(K_1^T,K_2^T,K_3^T,K_4^T,K_5^T)^T.
\end{array}
\]
Then the Jacobian of \(\vect f\mapsto\mathcal M(\mat L\vect f)\) is
\begin{equation}
    \mat E=\mat K\mat L.
\label{eq:true_equilibrium_jacobian}
\end{equation}
Moreover, \(\mat E^2=\mat E\),
\(\operatorname{rank}\mat E=3\), and
\(\ker\mat E=\ker\mat L\).
\end{lemma}

\begin{proof}
Differentiating the five equilibrium populations and then evaluating at
\(w_*\) gives the stated block rows.  Conservation implies
\(\mat L\mat K=I_3\), and hence
\(\mat E^2=\mat K(\mat L\mat K)\mat L=\mat E\).  Furthermore,
\(\operatorname{range}\mat E\subseteq\operatorname{range}\mat K\), while
\(\mat E\mat K=\mat K\) gives the reverse inclusion.  Thus
\(\operatorname{range}\mat E=\operatorname{range}\mat K\), and this space has
dimension three because \(\mat L\mat K=I_3\).  Finally,
\(\mat E\vect x=0\) implies
\(\mat L\vect x=\mat L\mat E\vect x=0\); the reverse inclusion follows from
\(\mat E=\mat K\mat L\).  Therefore \(\ker\mat E=\ker\mat L\).
\end{proof}

For each \(j=1,\ldots,5\),
\[
 \frac{\partial(M_1)_0}{\partial(f_j)_1}=\frac{\alpha}{2}.
\]
Thus \(\mat E\) contains density--momentum coupling blocks and differs from
\((a,a,a,a,1-4a)^T\vect 1_5^T\otimes I_3\).
The outer product over the velocity index has rank one, whereas the full
tensor product has rank three.

\subsection{Raw and modified moment variables}

Let \(\mat T=\Mtrans\otimes I_3\) and
\(\widehat{\mat S}=\mat S\otimes I_3\).  The equilibrium moments are
\[
    \vect m^{\rm eq}(w)
    =\left(w,\alpha A_1(w),\alpha A_2(w),0,0\right)^T.
\]
Consequently,
\begin{equation}
 \widehat{\mat E}:=\mat T\mat E\mat T^{-1}
 =\begin{pmatrix}
 I_3&0&0&0&0\\
 \alpha B_1&0&0&0&0\\
 \alpha B_2&0&0&0&0\\
 0&0&0&0&0\\
 0&0&0&0&0
 \end{pmatrix}.
\label{eq:equilibrium_jacobian_moments}
\end{equation}
For the collision increment
\(\vect Q(\vect f)=-\mat T^{-1}\widehat{\mat S}\mat T
(\vect f-\mathcal M(\mat L\vect f))\), its Jacobian at rest satisfies
\begin{equation}
 \widehat{\mat J}:=\mat T D\vect Q(\vect f_*)\mat T^{-1}
 =\begin{pmatrix}
 0&0&0&0&0\\
 s_1\alpha B_1&-s_1I_3&0&0&0\\
 s_2\alpha B_2&0&-s_2I_3&0&0\\
 0&0&0&-s_3I_3&0\\
 0&0&0&0&-s_4I_3
 \end{pmatrix}.
\label{eq:true_collision_jacobian_moments}
\end{equation}
The raw moment matrix therefore does not diagonalize the collision Jacobian.

To remove the blocks coupling the conserved variables to the flux moments,
introduce the lower triangular change of variables
\begin{equation}
 \mat R=
 \begin{pmatrix}
 I_3&0&0&0&0\\
 -\alpha B_1&I_3&0&0&0\\
 -\alpha B_2&0&I_3&0&0\\
 0&0&0&I_3&0\\
 0&0&0&0&I_3
 \end{pmatrix}.
\label{eq:adapted_transform}
\end{equation}
For a perturbation \(\delta\vect m\) about
\(\vect m_*=\mat T\vect f_*\), the modified nonconserved variables are
\(\delta m_1-\alpha B_1\delta m_0\),
\(\delta m_2-\alpha B_2\delta m_0\), \(\delta m_3\), and \(\delta m_4\).
They are the linearized nonequilibrium parts of the raw moments at \(w_*\).
Direct multiplication gives
\begin{equation}
 \mat R\widehat{\mat J}\mat R^{-1}
 =\diag(0,-s_1I_3,-s_2I_3,-s_3I_3,-s_4I_3).
\label{eq:adapted_diagonalization}
\end{equation}
Equation~\eqref{eq:adapted_diagonalization} shows that \(\mat R\mat T\)
diagonalizes the collision Jacobian, whereas \(\mat T\) alone does not.

\subsection{Exact nonlinear contraction of the kinetic residual}

Recall
\[
 \mat W=\diag(a,a,a,a,1-4a),\qquad 0<a<\frac14.
\]
For this moment basis, weighted orthogonality takes the following form.

\begin{proposition}[Weighted moment orthogonality]
\label{prop:orthogonality}
The matrix \(\mat D=\Mtrans\mat W\Mtrans^T\) is diagonal and positive
definite,
with
\begin{equation}
 \mat D=\diag\left(1,2a,2a,4a,\frac{4a}{1-4a}\right).
\label{eq:D2_weighted_moment_matrix}
\end{equation}
\end{proposition}

\begin{proof}
Opposite-direction odd rows are orthogonal to the even rows, and the two
flux rows have disjoint support.  The stress row is orthogonal to the
conserved and energy rows by cancellation.  Finally,
\(4a-(1-4a)4a/(1-4a)=0\) gives orthogonality between the conserved and energy
rows.  The diagonal entries follow by direct multiplication.
\end{proof}

For a state with positive density, define the kinetic residual
\begin{equation}
    \vect r(\vect f)=\vect f-\mathcal M(\mat L\vect f).
\label{eq:kinetic_residual}
\end{equation}
Although \(\vect r\) is nonlinear through \(\mathcal M\), conservation gives
\(\mat L\vect r(\vect f)=0\).

\begin{theorem}[Exact kinetic contraction of the MRT collision]
\label{thm:wl2_main}
Let \(0<a<1/4\) and \(s_k\in[0,2]\), \(k=1,\ldots,4\).  If
\(\mathcal C(\vect f)=\vect f+\vect Q(\vect f)\) denotes the local collision
map, then
\begin{equation}
 \|\vect r(\mathcal C(\vect f))\|_{\mat W^{-1}\otimes I_3}^2
 \le
 \|\vect r(\vect f)\|_{\mat W^{-1}\otimes I_3}^2.
\label{eq:exact_kinetic_contraction}
\end{equation}
More precisely, with \(\vect\eta=\mat T\vect r(\vect f)\), one has
\(\eta_0=0\) and
\begin{align}
 \|\vect r(\vect f)\|_{\mat W^{-1}\otimes I_3}^2
 &=\sum_{k=1}^4D_{kk}^{-1}\|\eta_k\|_2^2,\\
 \|\vect r(\mathcal C(\vect f))\|_{\mat W^{-1}\otimes I_3}^2
 &=\sum_{k=1}^4(1-s_k)^2D_{kk}^{-1}\|\eta_k\|_2^2.
\label{eq:exact_kinetic_identity}
\end{align}
Hence the contraction factor is at most
\(\max_{1\le k\le4}|1-s_k|\), and it is strictly below one when every
\(s_k\in(0,2)\).
\end{theorem}

\begin{proof}
Collision preserves \(w=\mat L\vect f\), so the equilibrium before and after
collision is the same.  Therefore
\[
 \mat T\vect r(\mathcal C(\vect f))
 =(I_{15}-\widehat{\mat S})\mat T\vect r(\vect f).
\]
The congruence identity
\(\mat T^{-T}(\mat W^{-1}\otimes I_3)\mat T^{-1}
=\mat D^{-1}\otimes I_3\) gives~\eqref{eq:exact_kinetic_identity}, and
\((1-s_k)^2\le1\) proves the result.
\end{proof}

Theorem~\ref{thm:wl2_main} is an exact estimate for the nonlinear local
collision map.  The entries of \(\mat D\) depend on the normalization of the
moment rows.  They are therefore not physical eigenvalues and do not define
invariant measures of anisotropy under changes of coordinates.

\subsection{Symmetrizer for the linearized collision operator}

The diagonal metric \(\mat W^{-1}\otimes I_3\) in distribution space does not
symmetrize the full collision Jacobian.  Equivalently,
\(\widehat{\mat J}(\mat D\otimes I_3)\) is not symmetric: its $(1,0)$ block is
\(s_1\alpha B_1\), whereas its $(0,1)$ block vanishes.  The same obstruction
occurs for \(B_2\).

Let
\[
 \mat\Gamma=\diag(\Gamma_0,\Gamma_1,\Gamma_2,\Gamma_3,\Gamma_4),
 \qquad \Gamma_k=\Gamma_k^T>0,
\]
and define
\begin{equation}
 \mat H=\mat T^T\mat R^T\mat\Gamma\mat R\mat T.
\label{eq:adapted_symmetrizer}
\end{equation}

\begin{theorem}[Symmetrizer for the full linearized collision operator]
\label{thm:adapted_collision_stability}
Let \(\mat A_c=I_{15}+D\vect Q(\vect f_*)\).  For
\(s_k\in[0,2]\), \(k=1,\ldots,4\), the matrix~\eqref{eq:adapted_symmetrizer}
is symmetric positive definite and symmetrizes both the collision Jacobian
and its one-step map:
\begin{equation}
 [D\vect Q(\vect f_*)]^T\mat H
 =\mat H D\vect Q(\vect f_*),
 \qquad
 \mat A_c^T\mat H=\mat H\mat A_c.
 \label{eq:adapted_self_adjointness}
\end{equation}
It also satisfies
\begin{equation}
    \mat A_c^T\mat H\mat A_c\preceq\mat H.
\label{eq:adapted_collision_inequality}
\end{equation}
More explicitly, if \(\vect z=\mat R\mat T\vect g\), then
\begin{equation}
 \|\mat A_c\vect g\|_{\mat H}^2-\|\vect g\|_{\mat H}^2
 =-\sum_{k=1}^4s_k(2-s_k)\,z_k^T\Gamma_kz_k.
\label{eq:adapted_energy_identity}
\end{equation}
If \(s_k\in(0,2)\) for every \(k\), equality holds if and only if
\(z_k=0\) for \(k=1,\ldots,4\), equivalently, if \(\vect g\) belongs to the
tangent space of the equilibrium manifold at \(\vect f_*\).
\end{theorem}

\begin{proof}
Because \(\mat T\) and \(\mat R\) are invertible and \(\mat\Gamma\) is
positive definite, the congruence in Eq.~\eqref{eq:adapted_symmetrizer} is
positive definite.
Equation~\eqref{eq:adapted_diagonalization} gives
\[
 \mat R\mat T\mat A_c\mat T^{-1}\mat R^{-1}
 =\diag(I_3,(1-s_1)I_3,\ldots,(1-s_4)I_3).
\]
The transformed matrix has real scalar diagonal blocks and therefore
commutes with \(\mat\Gamma\), proving
Eq.~\eqref{eq:adapted_self_adjointness}.  Computing the change in the
squared \(\mat\Gamma\)-norm gives
Eq.~\eqref{eq:adapted_energy_identity} and the matrix inequality.
\end{proof}

\begin{remark}[Limitation to the collision step]
\label{rem:collision_only}
The construction~\eqref{eq:adapted_symmetrizer} generally produces off-diagonal
blocks coupling distinct velocities.  Periodic streaming has Fourier symbol
\[
 \mat U(\vect\theta)=\diag
 (\e^{-\mathrm{i}\vect\theta\cdot e_1},\ldots,
  \e^{-\mathrm{i}\vect\theta\cdot e_5})\otimes I_3.
\]
A metric independent of the wave number makes \(\mat U(\vect\theta)\) an isometry for
every \(\vect\theta\) only if its blocks coupling distinct velocities vanish;
indeed, for a block matrix \(\mat H=(H_{ij})\),
\[
 [\mat U(\vect\theta)^*\mat H\mat U(\vect\theta)]_{ij}
 =\e^{\mathrm{i}\vect\theta\cdot(e_i-e_j)}H_{ij}.
\]
Theorem~\ref{thm:adapted_collision_stability} does not establish this
block-diagonal property.  Hence the collision inequality alone does not imply
stability of the coupled collision--streaming discretization.
Section~\ref{sec:algorithm} examines the corresponding full Fourier
amplification matrix.
\end{remark}

\section{D3N7 formulation and consistency}
\label{sec:d3n7}

The D3N7 lattice uses
\begin{align*}
 e_1&=-e_4=(1,0,0)^T,\qquad e_2=-e_5=(0,1,0)^T,\\
 e_3&=-e_6=(0,0,1)^T,\qquad e_7=0.
\end{align*}
Each population belongs to \(\mathbb R^4\), and
\(w=(\rho,q_1,q_2,q_3)^T\) with \(q_k=h\rho u_k\).
Throughout this section, \(0<a<1/6\), so all direction weights below are
positive.

\subsection{Equilibrium and moments}

For \(k=1,2,3\), the fluxes in conserved variables are
\begin{equation}
 (A_k)_0=q_k,\qquad
 (A_k)_j=\frac{q_kq_j}{\rho}+\delta_{kj}(\rho-\bar\rho),
 \quad j=1,2,3.
\label{eq:d3_flux_conserved}
\end{equation}
The equilibrium populations are
\begin{equation}
 M_k=aw+\frac{\alpha}{2}A_k,\qquad
 M_{k+3}=aw-\frac{\alpha}{2}A_k,\qquad
 M_7=(1-6a)w.
\label{eq:d3_equilibrium}
\end{equation}
With \(\gamma=6a/(1-6a)\), the moment matrix is
\begin{equation}
 \Mtrans=\begin{pmatrix}
 1&1&1&1&1&1&1\\
 1&0&0&-1&0&0&0\\
 0&1&0&0&-1&0&0\\
 0&0&1&0&0&-1&0\\
 1&-1&0&1&-1&0&0\\
 1&1&-2&1&1&-2&0\\
 1&1&1&1&1&1&-\gamma
 \end{pmatrix}.
\label{eq:M_trans_d3n7_rebuilt}
\end{equation}
Moment \(m_0\) is conserved; \(m_1,m_2,m_3\) are flux moments with
equilibrium values \(\alpha A_k\), \(k=1,2,3\).  The moments \(m_4,m_5\)
are trace-free stress combinations, and \(m_6\) is the energy moment.
The last three equilibrium moments vanish.

For
\(\mat W_3=\diag(a,a,a,a,a,a,1-6a)\), direct multiplication gives
\begin{equation}
 \Mtrans\mat W_3\Mtrans^T
 =\diag\left(1,2a,2a,2a,4a,12a,\frac{6a}{1-6a}\right).
\label{eq:D_d3n7_rebuilt}
\end{equation}
As in D2N5, these coefficients depend on the moment normalization and are
not physical eigenvalues.  Rescaling a moment row changes the corresponding
coefficient but leaves the collision operator unchanged.

\subsection{Relaxation rates and viscosity}

To preserve symmetry under permutations of the coordinate axes, we assign a common relaxation rate to
the three flux moments and a common rate to the two trace-free stress moments:
\begin{equation}
 \mat S=\diag(0,s_\nu,s_\nu,s_\nu,s_\sigma,s_\sigma,s_\varepsilon).
\label{eq:d3_independent_mrt_rates}
\end{equation}
The common flux rate yields an isotropic coefficient in the leading viscous
term.  The common stress rate preserves permutation symmetry on the
two-dimensional space of diagonal trace-free moments because a scalar rate
commutes with its representation; this rate does not enter the viscosity
relation.
The energy rate \(s_\varepsilon\) is independent of \(s_\sigma\).
Setting \(s_\sigma=s_\varepsilon\) recovers TRT, whereas unequal rates
give an MRT collision operator outside the TRT subfamily.

The viscosity relation is
\begin{equation}
 \nu=\frac{2a}{\alpha}\left(\frac1{s_\nu}-\frac12\right),
 \qquad
 s_\nu=\frac{2a}{\alpha\nu+a}.
\label{eq:d3_viscosity_rebuilt}
\end{equation}
For the equal weights \(a=1/7\), this becomes
\(s_\nu=2/(7\alpha\nu+1)\).

\begin{corollary}[Formal second-order consistency of D3N7]
\label{cor:d3_formal_consistency}
Let \(0<a<1/6\), let \(\alpha>0\) be fixed, and set
\[
 h=\Delta x,\qquad \Delta t=\alpha h^2,\qquad
 w=(\rho,h\rho u_1,h\rho u_2,h\rho u_3)^T,\qquad
 \rho=\bar\rho+h^2P .
\]
Assume that \(P\), \(\vect u\), and all derivatives required below are
uniformly bounded on an \(h\)-independent time interval, and that the D3N7
streaming Taylor expansion through parabolic degree four has a componentwise
\(O(h^5)\) remainder.  Let the relaxation rates be fixed as \(h\to0\) and
satisfy
\[
 0<s_{\min}\le s_\nu<2,
 \qquad
 s_{\min}\le s_\sigma,s_\varepsilon\le2.
\]
Assume further that the nonconserved Maxwell equation admits the expansion
\[
 \vect m_{\rm D3}(\lambda)
  =\vect m_{\rm D3}^{\rm eq}
   +\sum_{r=1}^{4}\lambda^rh^r\vect m_{\rm D3}^{(r)}
   +\vect{\mathcal R}_{m,{\rm D3}}^{(5)}(\lambda;h),
\]
with zero conserved corrections, uniformly bounded coefficients and required
derivatives, and a componentwise \(O(h^5)\) remainder uniformly for
\(0\le\lambda\le1\).  The conserved balance is imposed after setting
\(\lambda=1\).

Then the smooth, unforced interior D3N7 equivalent equations satisfy,
componentwise and with remainders uniform on the stated time interval,
\begin{align}
 \partial_t\rho+\nabla\!\cdot(\rho\vect u)&=O(h^2),\\
 \partial_t(\rho\vect u)
 +\nabla\!\cdot(\rho\vect u\otimes\vect u)
 +\nabla P-\nu\Delta(\rho\vect u)&=O(h^2),
\end{align}
where
\[
 \nu=\frac{2a}{\alpha}
       \left(\frac1{s_\nu}-\frac12\right).
\]
The result holds independently of the fixed values of
\(s_\sigma\) and \(s_\varepsilon\) in the stated interval.
\end{corollary}

\begin{proof}
The recursion in Eq.~\eqref{eq:cons_maxwell_recursion} and the conserved
projection in Eq.~\eqref{eq:cons_physical_projection} apply with the D3N7
operators, \(\vect e_0\in\mathbb R^7\), and \(I_4\) on each physical block.
Let \(\widehat{\mat D}_3\) denote the D3 counterpart of
Eq.~\eqref{eq:cons_Dhat_definition}.  Direct multiplication gives
\begin{align*}
 \vect m_{\rm D3}^{\rm eq}
   &=(w,\alpha A_1,\alpha A_2,\alpha A_3,0,0,0)^T,\\
 (\widehat{\mat D}_3\vect m_{\rm D3}^{\rm eq})_k
   &=2a\,\partial_kw\quad(k=1,2,3),\qquad
 (\widehat{\mat D}_3^{\,2}\vect m_{\rm D3}^{\rm eq})_0=2a\Delta w.
\end{align*}
The first two conserved coefficients are consequently
\[
 \vect c_{\rm D3}^{(1)}=\alpha\sum_{k=1}^3\partial_kA_k,
 \qquad
 \vect c_{\rm D3}^{(2)}=\alpha(\partial_tw-\nu\Delta w),
\]
with \(\nu\) from Eq.~\eqref{eq:d3_viscosity_rebuilt}.
With
\[
 \mat J_3=\diag(1,-1,-1,-1,1,1,1)\otimes I_4
\]
and
\[
 \widehat{\mat S}_3^\dagger
 =\diag(0,s_\nu^{-1},s_\nu^{-1},s_\nu^{-1},
        s_\sigma^{-1},s_\sigma^{-1},s_\varepsilon^{-1})\otimes I_4,
\]
one has
\(\mat J_3\widehat{\mat D}_3=-\widehat{\mat D}_3\mat J_3\) and
\([\mat J_3,\widehat{\mat S}_3^\dagger]=0\), even when
\(s_\sigma\ne s_\varepsilon\).  Set
\[
 \mat\Pi_{3,\tau}=\frac12(I_{28}+\tau\mat J_3),
 \qquad \tau\in\{+1,-1\}.
\]
Applying the induction in
Eq.~\eqref{eq:cons_projected_parity_induction} with the D3 operators shows
that the conserved coefficient of odd order is generated only by
\(\mat\Pi_{3,-}\vect m_{\rm D3}^{\rm eq}\), whereas that of even order is
generated only by \(\mat\Pi_{3,+}\vect m_{\rm D3}^{\rm eq}\).
Under the weakly compressible scaling, the
four physical components of the odd equilibrium part have orders
\((O(h),O(h^2),O(h^2),O(h^2))\); those of the nonconstant even part have
orders \((O(h^2),O(h),O(h),O(h))\).  Therefore,
\[
 \vect c_{\rm D3}^{(3)}=(O(h),O(h^2),O(h^2),O(h^2))^T,
 \qquad
 \vect c_{\rm D3}^{(4)}=(O(h^2),O(h),O(h),O(h))^T.
\]
After multiplication by their parabolic factors \(h^3\) and \(h^4\), and
collecting all contributions of parabolic degree at least five, the conserved
remainder is
\[
 \vect{\mathcal R}_{c,3}
   =(O(h^4),O(h^5),O(h^5),O(h^5))^T.
\]
At \(\lambda=1\), the D3 conserved balance is therefore
\[
 0=\alpha h\sum_{k=1}^3\partial_kA_k
   +\alpha h^2(\partial_tw-\nu\Delta w)
   +\vect{\mathcal R}_{c,3}.
\]
Its density and momentum components carry the physical prefactors
\(\alpha h^2\) and \(\alpha h^3\), respectively.  Division by these factors,
followed by \(\Delta\rho=h^2\Delta P\), gives the stated \(O(h^2)\) residuals.
\end{proof}

\subsection{D3N7 collision identities}

Let
\[
 \mat T_3=\Mtrans\otimes I_4,\qquad
 \mat L_3=\vect 1_7^T\otimes I_4,\qquad
 \widehat{\mat S}_3
 =\diag(0,s_\nu,s_\nu,s_\nu,s_\sigma,s_\sigma,s_\varepsilon)\otimes I_4,
\]
and denote the diagonal matrix on the right-hand side of
Eq.~\eqref{eq:D_d3n7_rebuilt} by \(\mat D_3\).  At the uniform rest state
\(w_*=(\bar\rho,0,0,0)^T\), set \(\vect f_*=\mathcal M(w_*)\) and
\(B_k=D_wA_k(w_*)\).  Define the block
lower triangular matrix \(\mat R_3\) by
\[
 (\mat R_3)_{jj}=I_4\quad(j=0,\ldots,6),\qquad
 (\mat R_3)_{k0}=-\alpha B_k\quad(k=1,2,3),
\]
with all remaining off-diagonal blocks equal to zero.

\begin{corollary}[D3N7 collision structure]
\label{thm:wl2_d3n7}
At rest, the D3N7 collision Jacobian satisfies
\[
 \mat R_3\mat T_3[D_{\vect f}\vect Q(\vect f_*)]
 \mat T_3^{-1}\mat R_3^{-1}
 =\diag(0,-s_\nu I_4,-s_\nu I_4,-s_\nu I_4,
          -s_\sigma I_4,-s_\sigma I_4,-s_\varepsilon I_4).
\]
For the kinetic residual
\(\vect r_3(\vect f)=\vect f-\mathcal M(\mat L_3\vect f)\) and
\(\vect\eta=\mat T_3\vect r_3(\vect f)\),
\begin{align*}
 \|\vect r_3(\vect f)\|_{\mat W_3^{-1}\otimes I_4}^2
 &=\sum_{k=1}^{6}(\mat D_3)_{kk}^{-1}\|\eta_k\|_2^2,\\
 \|\vect r_3(\mathcal C(\vect f))\|_{\mat W_3^{-1}\otimes I_4}^2
 &=\sum_{k=1}^{6}(1-s_k)^2(\mat D_3)_{kk}^{-1}\|\eta_k\|_2^2.
\end{align*}
Here \(s_{1,2,3}=s_\nu\), \(s_{4,5}=s_\sigma\), and
\(s_6=s_\varepsilon\).  Hence local collision is nonexpansive on the kinetic
residual when all nonconserved rates belong to \([0,2]\).
\end{corollary}

\begin{proof}
In moment coordinates, the D3N7 equilibrium Jacobian has only its conserved
block column nonzero: its first block is \(I_4\), and its next three blocks
are \(\alpha B_1,\alpha B_2,\alpha B_3\).  The transformation \(\mat R_3\)
eliminates these flux blocks and gives the stated diagonalization.  Collision
preserves \(\mat L_3\vect f\), so its equilibrium is unchanged during the
local update.  Equation~\eqref{eq:D_d3n7_rebuilt} then converts the diagonal
moment relaxation into the displayed weighted norm identities.
\end{proof}

The corresponding conditional identity for the axis-aligned
\(D_dN_{2d+1}\) velocity set in arbitrary dimension is stated and proved in
Section~S1 of the Supplementary Material.

The D3N7 velocity set contains only the six axis-aligned velocities and the rest
velocity.  Its mixed fourth-order velocity moments
\(\sum_{i=1}^7(\mat W_3)_{ii}(e_i)_p^2(e_i)_q^2\) therefore vanish for
\(p\ne q\), so this stencil lacks the fourth-order isotropy available
in the standard D3Q27 quadrature.  On that larger stencil,
Suga \emph{et al.}~\cite{Suga2015} develop an MRT model for turbulent flows,
whereas Geier \emph{et al.}~\cite{Geier2015} formulate and validate a cumulant
collision operator and compare it with moment-based MRT models.
The numerical study therefore samples both axis-aligned and diagonal Fourier wave
vectors and includes a three-dimensional exact solution.  Multidimensional
convergence is assessed directly rather than inferred from axis-aligned tests
alone.

\section{Selection of relaxation parameters by Fourier analysis}
\label{sec:algorithm}

The condition \(0\le s_k\le2\) yields the collision estimate in
Theorem~\ref{thm:wl2_main}, but it does not imply stability of the fully
discrete operator.  This section defines the Fourier-selection procedure
using the spectrum of the linearized collision--streaming map.  It combines
a finite grid of wave vectors with local maximization and uses no flow-error
data.  The separate flow calibrations are described in
Section~\ref{sec:numerical}.
For the periodic linearized problem, a perturbation proportional to
\(\exp(\mathrm{i}\vect\theta\cdot\vect x/h)\) retains its spatial form after
one time step, while its coefficient is multiplied by the amplification
matrix defined below.  The eigenvalues of this matrix describe the decay
and phase of the numerical modes.

\subsection{Amplification matrix}

For \(d\in\{2,3\}\), set \(n_c=d+1\), \(N=2d+1\), and
\(w_*=(\bar\rho,0,\ldots,0)^T\).  Define
\begin{equation}
 \mat L=\vect 1_N^T\otimes I_{n_c},\qquad
 \mat T=\mat M_d\otimes I_{n_c},\qquad
 \widehat{\mat S}=\mat S\otimes I_{n_c},
\label{eq:fourier_dimensioned_maps}
\end{equation}
where \(\mat M_2\) is the D2N5 matrix in Eq.~\eqref{eq:M_trans} and
\(\mat M_3\) is the D3N7 matrix in
Eq.~\eqref{eq:M_trans_d3n7_rebuilt}.  With
\(B_k=D_wA_k(w_*)\), define the block column
\(\mat K\in\mathbb R^{Nn_c\times n_c}\)
by
\begin{equation}
 K_k=aI_{n_c}+\frac{\alpha}{2}B_k,\quad
 K_{k+d}=aI_{n_c}-\frac{\alpha}{2}B_k\quad(k=1,\ldots,d),
 \qquad K_N=(1-2da)I_{n_c}.
\label{eq:fourier_equilibrium_tangent}
\end{equation}
Then \(\mat E=\mat K\mat L\) is the full equilibrium Jacobian in population
coordinates in either dimension, consistent with the D2N5 and D3N7 collision
identities derived above.  The derivative of the
post-collision map is
\begin{equation}
 \mat C(\mat S)=I-\mat T^{-1}\widehat{\mat S}\mat T(I-\mat E).
\label{eq:full_collision_matrix_fourier}
\end{equation}
For a nondimensional wave vector \(\vect\theta\in[-\pi,\pi]^d\), streaming is
represented by
\begin{equation}
 \mat U(\vect\theta)=\diag
 \left(\e^{-\mathrm{i}\vect\theta\cdot e_1},\ldots,
       \e^{-\mathrm{i}\vect\theta\cdot e_N}\right)\otimes I_{n_c}.
\label{eq:streaming_symbol}
\end{equation}
The one-step amplification matrix is therefore
\begin{equation}
 \mat G(\vect\theta;\mat S)=\mat U(\vect\theta)\mat C(\mat S).
\label{eq:full_amplification_matrix}
\end{equation}
Equation~\eqref{eq:full_amplification_matrix} retains the pressure and
momentum derivatives in \(\mat E\); replacing \(\mat E\) by a scalar velocity
projection changes the symbol and may change both the computed spectral radius
and the set of rates that satisfy the spectral criterion.

The viscosity relation fixes
\begin{equation}
 s_\nu=\frac{2a}{\alpha\nu+a}\in(0,2),
\label{eq:fourier_fixed_viscosity_rate}
\end{equation}
and the relaxation matrices have the forms
\begin{align}
 \mat S_{\rm D2}&=\diag(0,s_\nu,s_\nu,s_\sigma,s_\varepsilon),
 \label{eq:d2_candidate_rates}\\
 \mat S_{\rm D3}&=\diag(0,s_\nu,s_\nu,s_\nu,
                         s_\sigma,s_\sigma,s_\varepsilon).
 \label{eq:d3_candidate_rates}
\end{align}
The viscosity fixes \(s_\nu\), whereas \(s_\sigma\) and \(s_\varepsilon\)
are varied independently.  When \(s_\sigma\ne s_\varepsilon\), the scheme is a
general MRT choice outside the TRT subfamily.  When
\(s_\sigma=s_\varepsilon\), it belongs to TRT, with SRT as the special case
\(s_\sigma=s_\varepsilon=s_\nu\).

\subsection{Spectral test on the sampled wave vectors}

For the sampled test we use
\begin{equation}
 \Theta_{\rm stab}=\{j\pi/4:j=0,\ldots,4\}^{d}.
 \label{eq:stability_wavevector_set}
\end{equation}
For this finite set, define
\begin{align}
 R_{\max,\mathrm{samp}}&=\max_{\vect\theta\in\Theta_{\rm stab}}
   \spr\bigl(\mat G(\vect\theta;\mat S)\bigr),
 \label{eq:sampled_spectral_radius}\\
 R_{\mathrm{hf,samp}}&=
 \max_{\substack{\vect\theta\in\Theta_{\rm stab}\\
                  \|\vect\theta\|_2\ge\pi/2}}
 \spr\bigl(\mat G(\vect\theta;\mat S)\bigr),
 \label{eq:high_frequency_radius}
\end{align}
where \(\spr\) denotes the spectral radius.  Every candidate in the computed
10\% accuracy band, together with each additional candidate needed to
enforce the damping constraint below, is subjected to an additional
maximization over
\([0,\pi]^d\cap\{\|\vect\theta\|_2\ge\pi/2\}\).  Powell searches over the full
domain are initialized at eight points selected by a deterministic
diversity criterion.  We also search one representative of each
orbit generated by signed coordinate permutations of the box faces, edges,
and corners, apply bounded
optimization on one-dimensional boundary segments, and parametrize the inner
boundary \(\|\vect\theta\|_2=\pi/2\) in angular coordinates.  Let
\(R_{{\rm hf},{\rm loc}}\) be the largest spectral radius evaluated in these
searches and set
\begin{equation}
 R_{\rm scr}=\max\{R_{\max,\mathrm{samp}},R_{{\rm hf},{\rm loc}}\}.
\label{eq:screened_spectral_radius}
\end{equation}
All sampled high-frequency points are included among the initial evaluations;
hence \(R_{{\rm hf},{\rm loc}}\ge R_{\mathrm{hf,samp}}\).
A candidate passes this spectral test only if every local optimization
terminates successfully and \(R_{\rm scr}\le1+\tau\), where
\(\tau=5\times10^{-10}\).

The axis-aligned velocity set, the grouping of relaxation rates, and the
Jacobian at rest are invariant under signed coordinate permutations.
Consequently, for each matrix \(\mat O\) representing such a permutation,
there is an orthogonal signed permutation \(\mat\Pi_O\) acting jointly on
population and component indices such that
\begin{equation}
 \mat G(\mat O\vect\theta;\mat S)
  =\mat\Pi_O\mat G(\vect\theta;\mat S)\mat\Pi_O^{-1}.
\label{eq:signed_permutation_similarity}
\end{equation}
It therefore suffices to search the positive quadrant in D2 and the positive
octant in D3.  Although this symmetry reduction is exact, finite sampling and
local optimization do not prove that the maximum over the continuous domain
has been attained.  At \(\vect\theta=0\), the amplification matrix has unit
eigenvalues on \(\operatorname{range}\mat E\), the tangent space of the
equilibrium manifold.  For the candidate rates considered
here, all kinetic eigenvalues have modulus less than one, and hence
\(\spr(\mat G(0))=1\).  The origin is therefore included in
\(R_{\max,\mathrm{samp}}\).  For the selected parameter pair, we repeat the
calculation on \(25^2\) points in D2 or \(13^3\) points in D3.  We then
locally maximize the unrestricted spectral radius, starting from the eight
grid points with the largest sampled values.  At the point with the largest
value found by these searches, and at the largest high-frequency point on
the dense grid, we compute
\(\max_{0\le n\le50}\|\mat G^n\|_2\) to measure transient growth over the
first 50 time steps.

\subsection{Continuation of the transverse invariant subspace}

A basis of the shear eigenspace is not uniquely determined when the
corresponding eigenvalue has multiplicity greater than one.  Moreover,
matching to the nearest target eigenvalue may select a kinetic branch.  We therefore
continue the \((d-1)\)-dimensional transverse invariant subspace along each
Fourier ray.

For a unit direction \(\vect n\), let \(\mat P_\perp\) have orthonormal columns
spanning \(\vect n^\perp\), and initialize the continuation with
\begin{equation}
 \mat Q_T^0(\vect n)=\operatorname{orth}\!\left[
 \mat K\begin{pmatrix}0\\\mat P_\perp\end{pmatrix}\right].
\label{eq:transverse_equilibrium_tangent}
\end{equation}
Starting at \(\kappa_0=10^{-3}\), the \((d-1)\)-dimensional invariant
subspace at each radius is chosen to maximize its overlap with the previously
tracked subspace.  For a candidate index set \(I\), let \(\mat Q_I\) be an
orthonormal basis for its invariant subspace and define
\begin{equation}
 \omega(I)=\sigma_{\min}\!\left((\mat Q_T^{\rm old})^*\mat Q_I\right).
 \label{eq:tracking_overlap}
\end{equation}
Let \(I_*\) and \(I_{(2)}\) denote the candidates with the largest and
second-largest values of \(\omega(I)\).  We compute
\begin{equation}
 m=\omega(I_*)-\omega(I_{(2)}),\qquad
 g=\min_{\substack{j\in I_*\\k\notin I_*}}|\lambda_j-\lambda_k|,
 \qquad
 r=\|\mat G\mat Q_{I_*}
       -\mat Q_{I_*}(\mat Q_{I_*}^*\mat G\mat Q_{I_*})\|_2.
 \label{eq:invariant_subspace_diagnostics}
\end{equation}
A continuation is rejected if
\(\omega(I_*)<\cos(15^\circ)\), \(m<0.02\), \(g<10^{-7}\), or
\(r>10^{-8}\).  For fixed candidate subspaces, the overlap and residual are
independent of the orthonormal bases used to represent them.  The candidate
sets are constructed from individual eigenvectors, however, so a rotation
within a repeated eigenspace can change the competing candidates and the
margin $m$.  These thresholds are numerical diagnostics, not a
basis-independent selection rule.  Section~S2.2 of the Supplementary Material specifies the
construction of the candidate subspaces and the rule for equal overlap values
and reports rotation tests for the four
D3 rate choices.  The continuum multiplier is not used for branch identification.

For the tracked shear multipliers, the continuum reference is
\begin{equation}
 \lambda_{\rm NS}(\vect\theta)=
 \exp(-\alpha\nu\|\vect\theta\|_2^2)
\label{eq:continuum_shear_target}
\end{equation}
and the maximum normalized error is
\begin{equation}
 E_T=\max_{\vect\theta\in\Theta_{\rm acc}}
 \max_{1\le j\le d-1}
 \frac{|\lambda_{T,j}(\vect\theta)-\lambda_{\rm NS}(\vect\theta)|}
 {\max\{1-\lambda_{\rm NS}(\vect\theta),10^{-10}\}}.
\label{eq:normalized_shear_error}
\end{equation}
The damping and phase components are evaluated separately.  The continuum target
corresponds to one time step \(\Delta t=\alpha h^2\) for the physical Fourier
wave vector \(\vect\theta/h\); it is not an exact multiplier of the discrete
scheme.  The continuation starts at \(\kappa_0=10^{-3}\); this initialization
point is not included in the accuracy objective.  The accuracy set is
\begin{equation}
 \Theta_{\rm acc}
 =\{\kappa\vect n:\kappa=0.02j,\ j=1,\ldots,24,\
                     \vect n\in\mathcal N_d\}.
 \label{eq:accuracy_wavevector_set}
\end{equation}
For D2,
\begin{equation}
 \mathcal N_2
 =\{(\cos(j\pi/32),\sin(j\pi/32)):j=0,\ldots,8\}.
 \label{eq:d2_accuracy_directions}
\end{equation}
For D3, \(\mathcal N_3=\{\vect q/\|\vect q\|_2:\vect q\in\mathcal Q_3\}\),
where
\begin{equation}
\begin{split}
 \mathcal Q_3=\{&
 (1,0,0),(4,1,0),(2,1,0),(4,3,0),(1,1,0),\\
 &(4,1,1),(3,1,1),(2,1,1),(3,2,1),(4,3,1),
 (2,2,1),(3,3,2),(1,1,1)\}.
\end{split}
 \label{eq:d3_accuracy_directions}
\end{equation}
Thus \(E_T\) is evaluated only for
\(0.02\le\|\vect\theta\|_2\le0.48\).  Because \(E_T\) is a maximum rather
than an angular average, the multiplicities of symmetry orbits do not affect the
objective.  Damping outside this band is measured separately by
\(R_{{\rm hf},{\rm loc}}\).

\subsection{Finite search over the relaxation rates}

The admissible box for the free rates is
\begin{equation}
 \mathcal P=[0.05,1.95]^2,
 \qquad
 \mathcal P_{\rm TRT}=\{(s,s):0.05\le s\le1.95\}\subset\mathcal P.
 \label{eq:trt_selection_set}
\end{equation}
The D2 search combines a \(61\times61\) tensor grid, 1001 points on the TRT
diagonal, and 401 points on each boundary edge; the corresponding D3 counts
are \(31\times31\), 401, and 201.  Each connected component formed by feasible
sampled points is refined locally around its minimizers of \(E_T\) and
\(R_{{\rm hf},{\rm samp}}\).  Candidates entering the discrete 10\% accuracy
band, together with those needed to test sensitivity to the damping criterion,
are then
subjected to the local high-frequency maximization.  The spectral test is repeated
after every rate refinement.  Let \(\mathcal D_F\) be the resulting finite
candidate set for family \(F\), and let
\(\mathcal A_F\subset\mathcal D_F\) contain the candidates for which all
calculations completed and the sampled spectral, continuation, and
local high-frequency tests passed.  The candidate grids and local refinement
procedure are described in Section~S2.3 of the Supplementary Material;
the Fourier maximization strategy is outlined in Section~S2.1.
Define
\begin{equation}
 E_{F,\min}=\min_{x\in\mathcal A_F}E_T(x),\qquad
 \mathcal B_F=\{x\in\mathcal A_F:
 E_T(x)\le1.05E_{F,\min}\}.
 \label{eq:accuracy_band_candidate_set}
\end{equation}
Set
\begin{equation}
 R_{F,\min}^{\rm hf}=\min_{x\in\mathcal B_F}R_{{\rm hf},{\rm loc}}(x),
 \qquad
 \mathcal C_F=\{x\in\mathcal B_F:
 R_{{\rm hf},{\rm loc}}(x)\le R_{F,\min}^{\rm hf}+10^{-12}\}.
 \label{eq:accuracy_band_representative}
\end{equation}
The reported representative \(x_F^*\) is the lexicographic minimum over
\(\mathcal C_F\) of
\[
 \left(E_T,R_{\rm scr},
 \max\{|1-s_\sigma|,|1-s_\varepsilon|\},s_\sigma,s_\varepsilon\right).
\]
The absolute tolerance specifies numerical ties, which are resolved by the
lexicographic rule above.  The 5\% band specifies the permitted increase in
\(E_T\) above its minimum; flow errors are not used to select \(x_F^*\).
Accuracy bands of 0\%, 2\%, 5\%, and
10\%, together with the absolute damping bounds, are reported as sensitivity
data.  The TRT candidate with minimum error is inserted into the MRT candidate set, so
\(E_{{\rm MRT},\min}\le E_{{\rm TRT},\min}\) on these finite sets.  This set
inclusion concerns the selection objective and does not imply smaller errors
for every flow.

A second TRT reference used in Section~\ref{sec:numerical} is fixed by the
 complementary relation analyzed in Ref.~\cite{Zhao2026TRTPower}.  In the
present notation,
\begin{equation}
 s^-=s_\nu,\qquad s^+=s_\sigma=s_\varepsilon=2-s_\nu,
 \qquad
 \left(\frac1{s^+}-\frac12\right)
 \left(\frac1{s^-}-\frac12\right)=\frac14.
 \label{eq:power_trt_relation}
\end{equation}
We call this the complementary TRT choice.  The corresponding theorem of
power stability applies to linearization about a uniform rest state, with
positive equilibrium blocks and reversible homogeneous transport; it does
not cover the forced problems or those with nonhomogeneous boundary data
considered below.  Unlike the Fourier-selected TRT rate, this rate is fixed by theory
and is not fitted to a
Fourier objective or a flow error.

The amplification matrices have dimensions \(15\times15\) for D2N5 and
\(28\times28\) for D3N7.  Each spectral evaluation computes all eigenvalues
of a dense matrix and therefore has \(O(m^3)\) cost for matrix dimension
\(m\).  The rates are computed once before the flow simulation; this cost is
independent of the number of time steps.  All matrix
power norms are Euclidean operator norms in
population space.

\section{Numerical experiments}\label{sec:numerical}

Unless stated otherwise, the periodic computations use $\alpha=0.2$,
$\nu=0.1$, $a=1/5$ for D2N5, $a=1/7$ for D3N7, $T=0.05$, and
$\Delta t=\alpha h^2$.

Standalone Python code for the numerical examples, Fourier parameter
selection, and flow-error calibration is provided in the supplementary code
package.

\subsection{Fourier-selected rates and spectral checks}

The Fourier-selected rates were determined before the flow simulations by the finite
parameter search in Section~\ref{sec:algorithm}; the TRT candidates were restricted to
Eq.~\eqref{eq:trt_selection_set}.  Every candidate in the union of the 2\%,
5\%, and 10\% accuracy sets underwent both the sampled spectral test
and the local high-frequency maximization described in
Section~\ref{sec:algorithm}.  Candidates were excluded if a local optimizer
did not report convergence or if any ambiguity threshold in branch continuation was
violated.  Table~\ref{tab:selected-rates} reports the 5\%
representatives together with the fixed complementary TRT and SRT
reference choices.
These four a priori rate choices are independent of the flow-error data.

\begin{table*}[t]
\centering\scriptsize
\caption{Fourier-selected rates and fixed reference rates, with their
spectral quantities.  F-MRT and F-TRT denote the Fourier-selected
representatives in the 5\% accuracy band; C-TRT denotes complementary TRT,
and SRT is the single-rate reference.  None of these rates is flow-calibrated.  The viscosity rate $s_\nu$ is fixed by
Eq.~\eqref{eq:fourier_fixed_viscosity_rate}.  Smaller $E_T$ means a more
accurate shear multiplier, whereas smaller $R_{{\rm hf},{\rm loc}}$ means
stronger damping of the high-frequency modes considered in the search.}
\label{tab:selected-rates}
\setlength{\tabcolsep}{2.5pt}
\begin{tabular}{llcccrr}
\toprule
case & rate choice & $s_\nu$ & $s_\sigma$ & $s_\varepsilon$ & $E_T$ &
$R_{{\rm hf},{\rm loc}}$\\
\midrule
D2 periodic & F-MRT      & 1.818182 & 0.261895 & 1.936022 & $3.66398\times10^{-3}$ & 0.9695829\\
            & F-TRT      & 1.818182 & 0.513614 & 0.513614 & $5.55137\times10^{-3}$ & 0.9769886\\
            & C-TRT & 1.818182 & 0.181818 & 0.181818 & $1.87921\times10^{-2}$ & 0.9619235\\
            & SRT               & 1.818182 & 1.818182 & 1.818182 & $1.52277\times10^{-2}$ & 0.9561603\\
\addlinespace
D3 periodic & F-MRT      & 1.754386 & 0.345885 & 1.940599 & $3.51588\times10^{-3}$ & 0.9721035\\
            & F-TRT      & 1.754386 & 0.429965 & 0.429965 & $4.71977\times10^{-3}$ & 0.9768801\\
            & C-TRT & 1.754386 & 0.245614 & 0.245614 & $1.87036\times10^{-2}$ & 0.9627031\\
            & SRT               & 1.754386 & 1.754386 & 1.754386 & $2.16693\times10^{-2}$ & 0.9576467\\
\addlinespace
D2 channel  & F-MRT      & 1.951220 & 0.072884 & $1.950000^\dagger$ & $3.86764\times10^{-3}$ & 0.9960549\\
            & F-TRT      & 1.951220 & 0.076586 & 0.076586 & $7.40701\times10^{-3}$ & 0.9982313\\
            & C-TRT & 1.951220 & 0.048780 & 0.048780 & $1.90287\times10^{-2}$ & 0.9903552\\
            & SRT               & 1.951220 & 1.951220 & 1.951220 & $1.53538\times10^{-2}$ & 0.9888675\\
\bottomrule
\end{tabular}
\begin{flushleft}\footnotesize
$^\dagger$ Upper boundary of the prescribed rate box; this representative is
conditional on $\mathcal P=[0.05,1.95]^2$.
\end{flushleft}
\end{table*}

In all three parameter sets in Table~\ref{tab:selected-rates}, the Fourier-selected MRT point outside the TRT subfamily has both
smaller $E_T$ and smaller $R_{{\rm hf},{\rm loc}}$ than the TRT point selected
by the same rule.  For the two periodic cases, the reductions in $E_T$ are
34.0\% and 25.5\%, respectively.  The SRT and complementary TRT points have
smaller high-frequency spectral radii but
larger $E_T$.  Thus, stronger high-frequency damping and smaller shear-mode
error need not occur for the same rate set.

The parameter maps examined in the analysis retained only
candidates for which all spectral tests completed successfully; spectral
radii computed on the sampled grid alone were omitted.  Varying the accuracy
band among 0\%, 2\%, 5\%, and 10\%, and the damping level among 0.985, 0.990,
and 0.995, gave the same qualitative
comparison in the periodic cases.  We then re-evaluated the six Fourier-selected MRT
and TRT representatives, one pair for each parameter case, on a $25^2$ grid of
wavenumbers in D2 or a $13^3$ grid in D3 and locally maximized the
spectral radius over the full Fourier domain from eight initial points for
each representative.  All 48 optimizations converged, and none produced a
spectral radius above $1+3.2\times10^{-15}$.  At the resulting critical and
high-frequency points, $\max_{0\le n\le50}\|\mat G^n\|_2$ ranged from 1.412
to 2.425.  Section~S2.1 of the Supplementary Material gives these additional
spectral checks and the unstable example below.

For the D2 boundary parameters with $\nu=0.025$, the pair of free rates
$(s_\sigma,s_\varepsilon)=(1,1)$ has sampled spectral radius 1.008141 and
locally refined radius $1.013358879>1$, attained near
$\vect\theta=(0,2.098488)$.  It is therefore rejected although its branch
continuation is unambiguous.  The Fourier-selected MRT rates for the boundary tests
are $(s_\nu,s_\sigma,s_\varepsilon)=(1.9512195122,0.072884114583,1.95)$,
whereas the Fourier-selected TRT point is interior to its one-dimensional search
interval.

\paragraph{Sensitivity to the rate bounds and relation to flow errors}
Extending the upper bound on the free rates from 1.95 to 1.99 changes the D2 periodic
representative only slightly: the changes in $s_\sigma$ and $s_\varepsilon$
are less than $3.0\times10^{-4}$ and $6.7\times10^{-3}$, respectively.  Within
the connected D3 component containing the minimum of $E_T$, increasing
$s_\varepsilon$ continues to reduce $E_T$ but increases the high-frequency spectral radius,
so the point selected from the 5\% accuracy band remains interior.  For the channel parameters, a
continuation over 401 candidates places both the point with minimum error and
the 5\% representative at the new
upper bound on the energy rate.  The minimum $E_T$, however, changes by less than
0.08\%.  The flow comparisons use the original rate box; the selected
channel rates are therefore conditional on the prescribed bounds.

We also evaluated a fixed set of 24 rate pairs on four plane waves from two
D3 direction orbits not related by signed coordinate permutations to the
original accuracy directions.  The values of $E_T$ on the original and additional direction sets
have Spearman correlation 0.966, whereas the correlation between $E_T$ on the
original direction set and the maximum $n_g=32$ velocity error over the four
additional waves is only 0.167.  A separate calculation with the complete
propagator $\mat G(\vect\theta)^{32}$, including equilibrium initialization
and macroscopic projection, reproduces, for each rate pair, the maximum
velocity error over the four waves to within $7.93\times10^{-8}$.
Thus, although the values of $E_T$ are strongly
correlated between the two direction sets, $E_T$ alone does not predict
finite-time flow errors, which also depend on modal projection and kinetic
transients.  The complete finite-time errors over the candidate set
constructed using the original $E_T$ criterion were no smaller for MRT than
for TRT.  The rate choices used in the flow tests remain unchanged.
The correlations and the complete-propagator comparison are reported in
Section~S4 of the Supplementary Material; the supplementary code implements
the candidate calculations and rate-bound tests.

\subsection{Periodic exact flows}

Both tests are periodic on \(\Omega_d=[0,1]^d\), with \(n_g\) cell-centered
points per coordinate, \(h=1/n_g\), and uniform grid \(\mathcal G_h\).  The D2
Taylor--Green vortex~\cite{Taylor1937} is
\begin{align}
u_x&=-U_0\cos(kx)\sin(ky)\e^{-2\nu k^2t},&
u_y&= U_0\sin(kx)\cos(ky)\e^{-2\nu k^2t},\\
P&=-\frac{U_0^2}{4}\{\cos(2kx)+\cos(2ky)\}\e^{-4\nu k^2t},
\end{align}
where $k=2\pi$ and $U_0=0.05$.  The three-dimensional
Arnold--Beltrami--Childress (ABC) test is defined by
\begin{align}
\vect b(x,y,z)&=\begin{pmatrix}
\sin(kz)+\cos(ky)\\ \sin(kx)+\cos(kz)\\ \sin(ky)+\cos(kx)
\end{pmatrix},\\
\vect u&=A(t)\vect b,\qquad A(t)=0.02\e^{-\nu k^2t},\\
P&=-\tfrac12A(t)^2\left(|\vect b|^2-3\right).
\end{align}
The identities $\nabla\cdot\vect b=0$,
$\nabla\times\vect b=k\vect b$, and
$\Delta\vect b=-k^2\vect b$ imply that the nonlinear term is
balanced by the pressure gradient and that the stated pair is an unforced
solution.

The exact pressure and the numerical pressure
\(P_h=(\rho_h-1)/h^2\) are shifted to zero mean, and density is initialized as
$\rho=1+h^2P$.  Table~\ref{tab:periodic-summary} gives the relative errors on
the finest grid and the orders fitted over four grids for SRT, Fourier-selected TRT,
complementary TRT, and the Fourier-selected MRT.  The
viscosity, flux rate, grid, time step, initial data, and final time are
identical within each comparison.  The supplementary code computes these
quantities at every grid level for the four rate choices listed in
Table~\ref{tab:periodic-summary}.
For a scalar or vector grid field \(v\), define the discrete root-mean-square
norm
\begin{equation}
 \|v\|_{2,h}=\left(
  \frac{1}{|\mathcal G_h|}\sum_{x\in\mathcal G_h}|v(x)|_2^2
 \right)^{1/2}.
 \label{eq:numerical_discrete_l2}
\end{equation}
The periodic errors are
\(E_u=\|\vect u_h-\vect u\|_{2,h}/\|\vect u\|_{2,h}\) and
\(E_P=\|P_h-P\|_{2,h}/\|P\|_{2,h}\).
We also report
\begin{equation}
 E_\rho=\|\rho_h-\rho\|_{2,h},\qquad
 E_{\rm div}=\frac{\|\nabla_h\!\cdot\vect u_h\|_{2,h}}
 {k\|\vect u\|_{2,h}},\qquad
 E_m=|\langle\rho_h\rangle_h-\langle\rho_h(0)\rangle_h|,
\label{eq:periodic_auxiliary_diagnostics}
\end{equation}
where \(\nabla_h\) is the centered periodic difference and
\(\langle\cdot\rangle_h\) is the grid average.

\begin{table*}[t]
\centering\footnotesize
\caption{Errors for periodic exact solutions and observed orders.  The finest grids
are $n_g=32$ in D2 and $n_g=24$ in D3.  The orders are least-squares slopes of
$\log E$ versus $\log h$ over all four grids.  F-MRT, F-TRT, C-TRT,
and SRT use the corresponding periodic rates in Table~\ref{tab:selected-rates}.}
\label{tab:periodic-summary}
\setlength{\tabcolsep}{8pt}
\begin{tabular}{@{}llrrrr@{}}
\toprule
Case & rate set & $E_u^{\rm fine}$ & $r_u$ & $E_P^{\rm fine}$ & $r_P$\\
\midrule
D2 & SRT
 & 0.0071422 & 1.9946 & 0.014762 & 2.0073\\
& F-TRT
 & 0.0069259 & 1.9895 & 0.014065 & 2.0052\\
& C-TRT
 & 0.0063828 & 1.9859 & 0.012366 & 1.9653\\
& F-MRT
 & 0.0071474 & 2.0008 & 0.014183 & 1.9905\\
\midrule
D3 & SRT
 & 0.0060869 & 2.0356 & 0.01232 & 2.0433\\
& F-TRT
 & 0.0048581 & 1.9701 & 0.0097741 & 1.9429\\
& C-TRT
 & 0.0036723 & 1.9440 & 0.0073203 & 1.9032\\
& F-MRT
 & 0.0045798 & 2.0546 & 0.0091671 & 2.0164\\
\bottomrule
\end{tabular}
\begin{flushleft}\footnotesize
D2 grids use $n_g=(12,16,24,32)$ and D3 grids use
$n_g=(8,12,16,24)$.  The supplementary code computes the errors on all grids.
\end{flushleft}
\end{table*}

\begin{figure*}[t]
\centering
\includegraphics[width=0.94\textwidth]{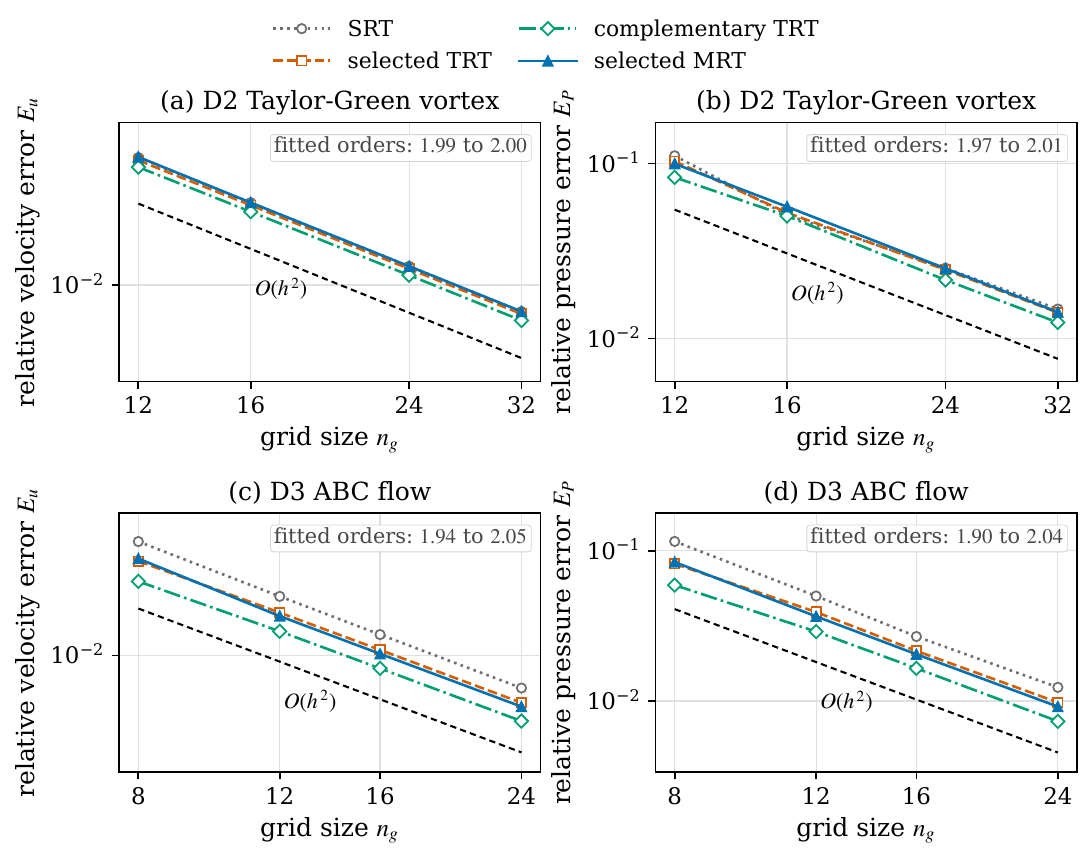}
\caption{Grid convergence for the D2 Taylor--Green and D3 ABC exact flows.
All collision families use the same viscosity, grid, time step, initial data,
and final time within each problem.  In the legends, ``selected MRT'' and
``selected TRT'' mean Fourier-selected rates (F-MRT and F-TRT).
Dashed black segments indicate
$O(h^2)$; the annotated values are least-squares orders over all four grids.}
\label{fig:core-second-order}
\end{figure*}

The least-squares orders in Table~\ref{tab:periodic-summary} range from 1.9032 to
2.0546.  Figure~\ref{fig:core-second-order} shows the velocity and pressure
convergence curves for the D2 Taylor--Green vortex in panels (a,b) and the
D3 ABC flow in panels (c,d).  In D3, the Fourier-selected MRT rates give smaller velocity and pressure
errors than the Fourier-selected TRT on the three finer grids; at
$n_g=24$, the reductions are 5.7\% and 6.2\%, respectively.  The
complementary TRT choice nevertheless has the smallest velocity and pressure
errors on every listed D2 and D3 grid.  Hence the ordering of the finite
Fourier objectives does not determine the error ordering for these nonlinear
exact solutions.

Section~S3 and Table~S1 of the Supplementary Material report the density
error, discrete divergence, and mass drift at every grid level for Fourier-selected MRT.

\subsection{Oblique transverse plane waves}

To test the order of the interior scheme in directions not used to evaluate
$E_T$ at small wave numbers,
let $\vect\ell\in\mathbb Z^3$, choose a constant unit vector $\vect p$
with $\vect\ell\cdot\vect p=0$, and set
\begin{equation}
 \vect u(\vect x,t)=A\vect p\sin(2\pi\vect\ell\cdot\vect x)
 \exp\{-\nu(2\pi)^2|\vect\ell|^2t\},
 \qquad P(\vect x,t)=0.
 \label{eq:oblique_shear_exact}
\end{equation}
This is an exact unforced Navier--Stokes solution because
$\nabla\cdot\vect u=0$, $(\vect u\cdot\nabla)\vect u=0$, and
$\Delta\vect u=-(2\pi)^2|\vect\ell|^2\vect u$.  We use
$\vect\ell=(3,1,0)$ with
$\vect p_1=(1,-3,0)/\sqrt{10}$ and $\vect p_2=(0,0,1)$, and independently
repeat the test for $\vect\ell=(4,2,1)$.  Neither direction belongs to
the orbit generated by signed coordinate permutations of any vector in
Eq.~\eqref{eq:d3_accuracy_directions}.  These directions occur in the separate
high-frequency spectral test but not in the accuracy objective $E_T$; no
finite-time plane-wave error was used to tune the relaxation rates.

\begin{figure*}[t]
\centering
\includegraphics[width=0.95\textwidth]{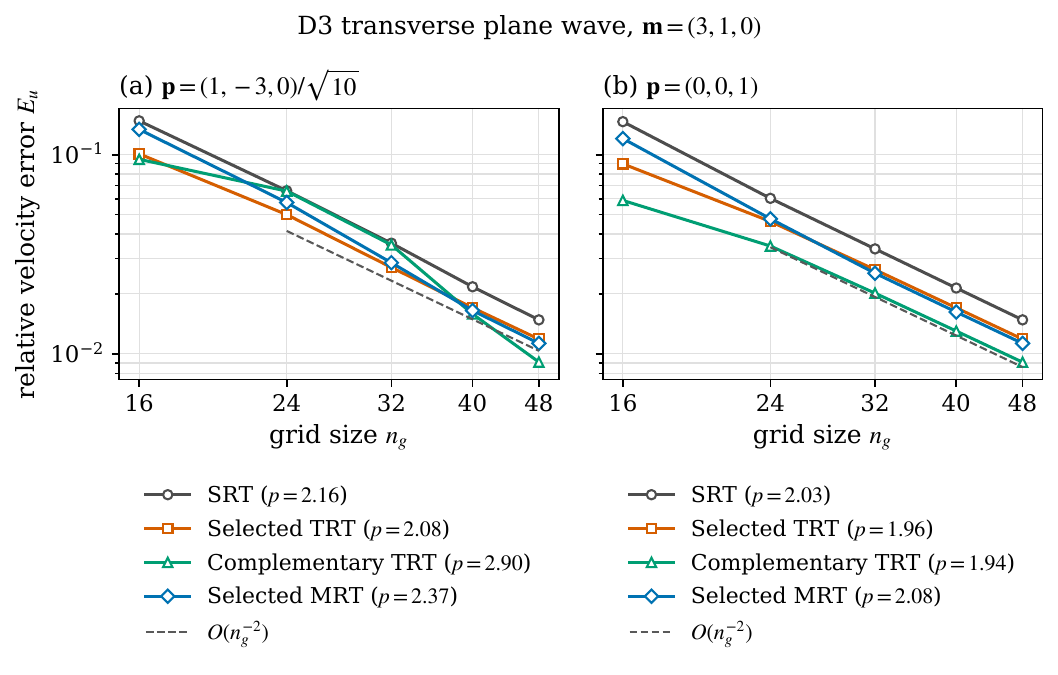}
\caption{Grid convergence for two polarizations of a D3 transverse plane wave
with $\vect\ell=(3,1,0)$.  The orbit of this wave vector is absent from the
direction set used for $E_T$, although it is represented in the
spectral test; all four rate choices were fixed before this test.
The legend labels ``Selected MRT'' and ``Selected TRT'' denote F-MRT and
F-TRT; the remaining choices are complementary TRT and SRT.
The displayed orders are least-squares fits over $n_g=24,32,40,48$; dashed
segments indicate $O(n_g^{-2})$.}
\label{fig:heldout-oblique-shear}
\end{figure*}

The computations use $A=0.02$, $T=0.00625$, and
$n_g=(16,24,32,40,48)$; the remaining parameters are those stated at the
beginning of this section.  The corresponding nondimensional wave numbers also extend
beyond the accuracy band $\|\vect\theta\|_2\le0.48$: for
$\vect\ell=(3,1,0)$, their magnitudes range from 1.242 to 0.414, whereas
$\vect\ell=(4,2,1)$ still gives 0.600 on the finest grid.  Thus these tests
examine additional directions and wave-number ranges together.
Figure~\ref{fig:heldout-oblique-shear} reports the
refinement results for all four rate choices with $\vect\ell=(3,1,0)$.

In Figure~\ref{fig:heldout-oblique-shear}, Fourier-selected MRT has 4.8\% smaller
error than Fourier-selected TRT at $n_g=48$ for both
polarizations, but complementary TRT has about 19.5\% smaller error than
Fourier-selected MRT.  Most slopes fitted over the finer grids are close to two.  We
regard the larger values for the first polarization of Fourier-selected MRT and
complementary TRT as pre-asymptotic and not as evidence of order greater
than two.  As a separate
check, the two slopes for Fourier-selected MRT with $\vect\ell=(4,2,1)$ are 1.966 and
1.962 over the same four finer grids.  Section~S4 and Table~S2 of the
Supplementary Material report the
refinement results for all four rate choices with $\vect\ell=(3,1,0)$ and the
Fourier-selected MRT results for $\vect\ell=(4,2,1)$; the supplementary code
implements these comparisons.  These tests provide further numerical evidence for the
second-order interior analysis.  The differences between the rate choices
concern error constants rather than convergence order.

\subsection{Nonlinear periodic flows with independent references}

The exact solutions above test the convergence order, but not general
nonlinear mode coupling.  We next consider two nonlinear flows and compare
them with independently computed reference solutions.  The relaxation rates
are the Fourier-selected and fixed reference choices in
Table~\ref{tab:selected-rates}, together with a separate TGV-calibrated pair
of MRT and TRT choices used only in the three-dimensional Taylor--Green tests.
The double-shear tests use only the four a priori D2 choices.

\paragraph{Three-dimensional Taylor--Green vortex}
On the unit three-torus, the initial velocity is
\begin{equation}
 \vect u_0=A\bigl(
 \sin(2\pi x)\cos(2\pi y)\cos(2\pi z),
 -\cos(2\pi x)\sin(2\pi y)\cos(2\pi z),0\bigr)^T.
 \label{eq:nonlinear_tgv3d_initial}
\end{equation}
We take $A=4$, $\nu=0.1$, and $T=0.1$.  The initial pressure $P_0$ is the
mean-zero spectral solution of
\begin{equation}
 \Delta P_0=-\nabla\!\cdot[(\vect u_0\!\cdot\nabla)\vect u_0],
 \qquad \rho_h(0)=1+h^2P_0.
 \label{eq:nonlinear_pressure_preparation}
\end{equation}
The independent reference uses a Fourier pseudospectral discretization with
the two-thirds dealiasing rule on a $48^3$ grid.  Time integration uses a
fourth-order Runge--Kutta method with an integrating factor and time step
$2\times10^{-4}$.  Comparisons of this reference with a $36^3$ calculation
and with a calculation using twice the time step give maximum relative differences of
$1.91\times10^{-10}$ in velocity and $4.89\times10^{-9}$ in vorticity.  The nonlinear reference differs
by 13.3\% from the viscously decayed initial Fourier modes.  For both
nonlinear tests, $E_u$ and $E_\omega$
denote the relative discrete $L^2$ errors of velocity and vorticity,
respectively.  The reference velocity series and its spectral curl are
evaluated directly at the VLBM cell centers on every comparison grid; VLBM
vorticity is obtained by Fourier differentiation on that grid.
Figure~\ref{fig:nonlinear-tgv3d}(a) shows a slice of the reference vorticity;
Section~S5 of the Supplementary Material gives the reference refinement checks.

In addition to the four a priori rate choices, TGV calibration minimizes
the final-time velocity and vorticity error objective
$J=\max(E_u,E_\omega)$ at $A=4$ on the single coarse grid $n_g=12$.  Refinement over 423 MRT candidates and 277
candidates restricted to the TRT diagonal gives
\((s_\sigma,s_\varepsilon)=(0.4953125,0.0796875)\) for MRT.  For TRT, the
corresponding pair is \((0.1750042,0.1750042)\).  These rates are then kept fixed
for the computations on $n_g=16,20,24$.  Both calibrated pairs also satisfy
the spectral and branch-continuation criteria defined above; the calibration
objective is $J$.
The comparison uses the minimizers of $J$ over the two finite candidate
sets; the two parameter searches have different computational costs.
Table~\ref{tab:nonlinear-tgv3d} reports the errors on the finest grid and
the fitted convergence orders for all six rate choices.
Figure~\ref{fig:nonlinear-tgv3d}(b) compares the two calibrated choices on
the finer grids, which were not used to fit the rates.  Complete refinement
data are given in Tables~S3 and S4 of the Supplementary Material.

\begin{table*}[t]
\centering\footnotesize
\caption{Nonlinear D3 Taylor--Green errors against the independent reference.
The reported errors are on $n_g=24$.  F-MRT and F-TRT are Fourier-selected;
C-TRT and SRT use fixed reference rates.  Their orders use
$n_g=12,16,20,24$.  TGV-calibrated rates minimize $\max(E_u,E_\omega)$ at
$A=4$, $n_g=12$; their orders use only $n_g=16,20,24$, which were excluded
from calibration.}
\label{tab:nonlinear-tgv3d}
\setlength{\tabcolsep}{4.5pt}
\begin{tabular}{@{}lccrrrr@{}}
\toprule
Rate choice & $s_\sigma$ & $s_\varepsilon$ & $E_u$ & $r_u$ & $E_\omega$ & $r_\omega$\\
\midrule
SRT & 1.754386 & 1.754386 & 0.010136 & 2.020 & 0.012715 & 1.988\\
F-TRT & 0.429965 & 0.429965 & 0.008356 & 2.022 & 0.010798 & 1.991\\
C-TRT & 0.245614 & 0.245614 & 0.006870 & 2.024 & 0.009389 & 2.020\\
F-MRT & 0.345885 & 1.940599 & 0.009704 & 2.010 & 0.012090 & 1.987\\
\addlinespace
TGV-calibrated TRT & 0.175004 & 0.175004 & 0.005815 & 1.995 & 0.008743 & 2.072\\
TGV-calibrated MRT & 0.495313 & 0.079688 & 0.003640 & 2.085 & 0.005800 & 2.120\\
\bottomrule
\end{tabular}
\end{table*}

\begin{figure*}[t]
\centering
\includegraphics[width=\textwidth]{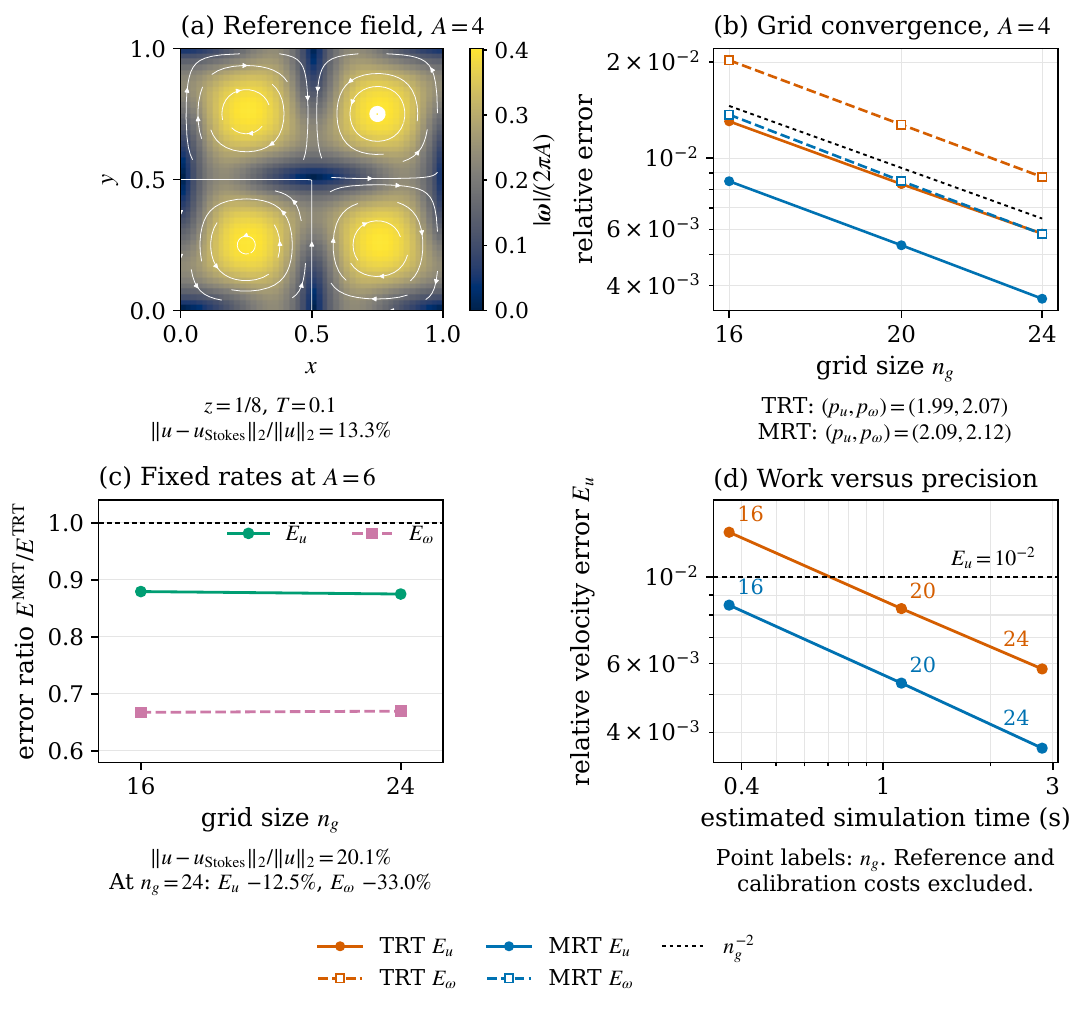}
\caption{D3 Taylor--Green vortex: reference flow, errors, and estimated
simulation time.  Panel (a) shows a reference vorticity slice at $A=4$.
The MRT and TRT labels in panels (b)--(d) refer exclusively to the
TGV-calibrated rates, not to F-MRT or F-TRT.
Panel (b) gives errors on grids not used to fit the rates, which were
obtained only at $n_g=12$.  Panel (c) uses the same rates at $A=6$, an
amplitude not used in the fit.  Panel (d) shows the estimated simulation
time required to attain a given accuracy on the finer grids at $A=4$;
the time excludes computation of the reference solution and calibration of
the relaxation rates.}
\label{fig:nonlinear-tgv3d}
\end{figure*}

Among the four a priori rate choices in Table~\ref{tab:nonlinear-tgv3d},
complementary TRT gives the smallest velocity and vorticity errors;
the Fourier-selected MRT rates do not improve on either TRT choice.
The coarse-grid fit gives a different comparison.  On the $n_g=12$ grid
used to fit the rates, the TGV-calibrated MRT value of $J$ is 28.9\% lower than
the TGV-calibrated TRT value (Section~S5 of the Supplementary Material).
On the three finer grids in Figure~\ref{fig:nonlinear-tgv3d}(b), both
calibrated choices give fitted orders close to two.  At $n_g=24$, TGV-calibrated MRT reduces the TGV-calibrated TRT velocity and vorticity errors by 37.4\% and
33.7\%, respectively, as shown in Table~\ref{tab:nonlinear-tgv3d}.

Figure~\ref{fig:nonlinear-tgv3d}(c) tests the same calibrated rates at
$A=6$, using a separate reference computed after the rates were fixed.
At $n_g=24$, TGV-calibrated MRT reduces the TGV-calibrated TRT velocity and
vorticity errors by 12.5\% and 33.0\%, respectively.  The reference departure
from linear viscous decay is 20.1\%; the reference checks and the errors on
all tested grids are reported in Section~S5 of the Supplementary Material.

Figure~\ref{fig:nonlinear-tgv3d}(d) compares the estimated simulation time
at $A=4$.  For the implementations used here, the median TGV-calibrated TRT and
MRT costs are both about 704 ns per lattice site and time step; Section~S5
of the Supplementary Material describes the timing measurements.
A velocity error below 1\% is first attained on the tested grids by TGV-calibrated MRT at $n_g=16$ and TGV-calibrated TRT at $n_g=20$, corresponding to estimated
simulation times of 0.369 and 1.126 s.  These times exclude the 53.0 s parameter search,
which includes 35.4 s of candidate solves for MRT and 17.0 s for TRT, and depend on the
implementation and grids used.  In this comparison, TGV-calibrated MRT requires less
simulation time because it attains the prescribed accuracy on a coarser grid.

\begin{figure*}[t]
\centering
\includegraphics[width=\textwidth]{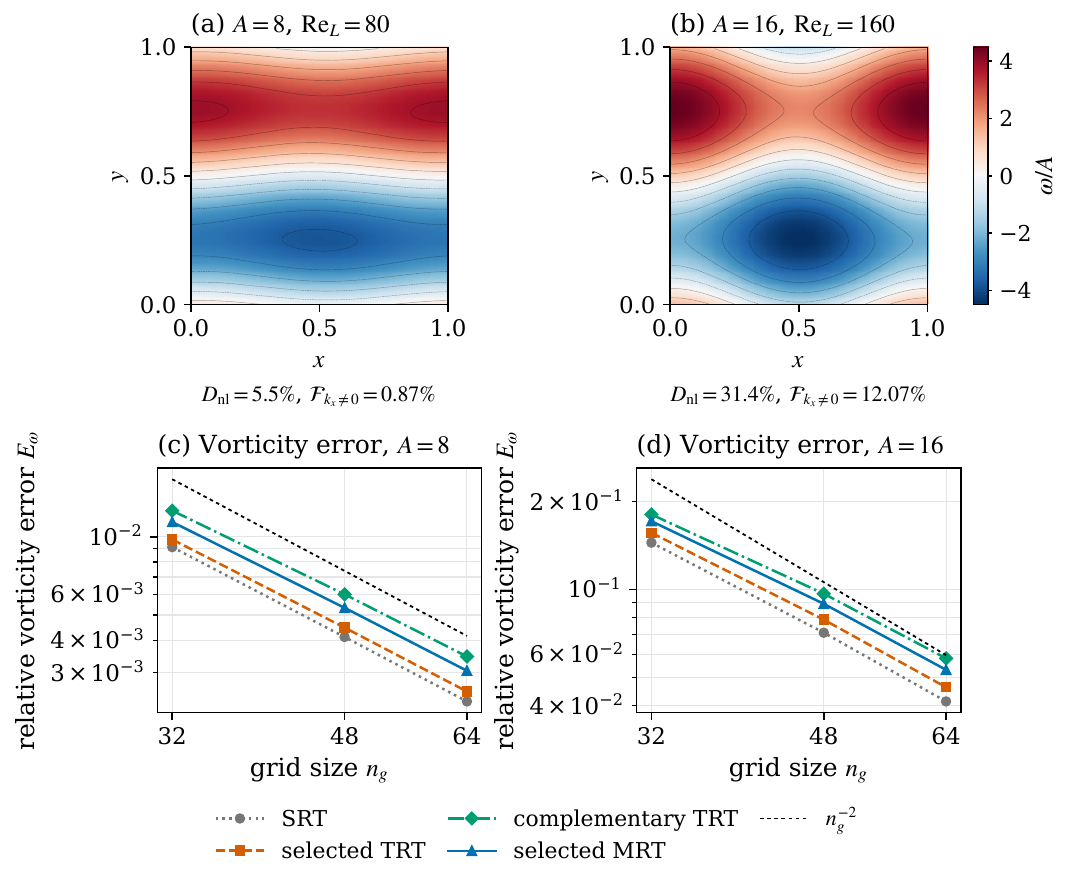}
\caption{Nonlinear roll-up of a double shear layer at $T=0.2$.  The upper panels show
normalized reference vorticity for $A=8$ and 16; their annotations give the
departure from purely diffusive evolution and the final fraction of enstrophy in
modes with nonzero streamwise wave number.  The lower panels show VLBM vorticity errors against
the independent reference solution.  The ``selected'' legend entries are
the Fourier-selected D2 rates; all four choices are a priori rates from
Table~\ref{tab:selected-rates}.  Full velocity and vorticity data are
reported in Section~S6 and Table~S5 of the Supplementary Material.}
\label{fig:double-shear-rollup}
\end{figure*}

\paragraph{Roll-up of a two-dimensional double shear layer}
To examine stronger nonlinear effects, we use the following initial condition for
a double shear layer
on the unit torus~\cite{BellColellaGlaz1989}:
\begin{equation}
 u_x(x,y)=
 \begin{cases}
 A\tanh[30(y-1/4)],&y\le1/2,\\
 A\tanh[30(3/4-y)],&y>1/2,
 \end{cases}
 \qquad u_y(x,y)=0.05A\sin(2\pi x).
 \label{eq:double_shear_initial_data}
\end{equation}
Using $A$ and the side length of the unit box as velocity and length scales gives
$\mathrm{Re}=A/\nu=80$ and 160 for $A=8$ and 16.  We use $T=0.2$,
$n_g=(32,48,64)$, and the four a priori D2 rate choices.  For each amplitude,
$P_0$ is the mean-zero spectral solution in
Eq.~\eqref{eq:nonlinear_pressure_preparation}, evaluated with the initial
velocity in Eq.~\eqref{eq:double_shear_initial_data}, and
$\rho_h(0)=1+h^2P_0$.  The dealiased reference, computed on a $192^2$ grid with
$\Delta t=5\times10^{-5}$, agrees with independent spatial and temporal checks to $8.92\times10^{-10}$
in both velocity and vorticity.
Figure~\ref{fig:double-shear-rollup}(a,b) shows the reference vorticity at
the two amplitudes, and panels (c,d) show the corresponding vorticity errors
under grid refinement.  Table~S5 of the Supplementary Material gives both
velocity and vorticity errors and their fitted orders.

For $A=8$ in Figure~\ref{fig:double-shear-rollup}(a,c), the reference
departure from purely diffusive evolution is 5.46\%,
and the fitted velocity orders in Table~S5 range from 1.86 to 1.97.  At $A=16$, nonlinear
mode transfer is much stronger (Figure~\ref{fig:double-shear-rollup}(b,d)):
the departure is 31.4\%, and the fraction of
enstrophy in modes with nonzero streamwise wave number increases from
0.062\% to 12.1\%.  The fitted orders 1.63--1.80 indicate that these grids are
not yet in a common second-order regime.  SRT gives the smallest velocity and
vorticity errors in Table~S5.  Section~S6 of the Supplementary Material also
notes that complementary TRT has the smallest energy and enstrophy errors
on the finest grid.  Thus the MRT rates selected by Fourier analysis
do not minimize the errors in this nonlinear roll-up.

\subsection{Dirichlet channel boundaries}

Relaxation rates can affect boundary errors in lattice Boltzmann
schemes~\cite{Ginzburg2005}.  We consider D2 channels with four
periodic cells in the streamwise direction and $n_g$ cells in the wall-normal
direction.  The initial state, body force, and wall data are streamwise
uniform, and the numerical update preserves this invariant subspace; hence the
extent of four streamwise cells does not affect these one-dimensional solutions.
These calculations use $a=1/5$, $\alpha=0.2$, $\nu=0.025$, and
$\Delta t=\alpha h^2$.  Collision is evaluated directly in moment space.  The
comparison includes the Fourier-selected MRT, the Fourier-selected TRT
rate in Table~\ref{tab:selected-rates}, and SRT with
$s_\sigma=s_\varepsilon=s_\nu$.  All three use the same forcing, wall update,
grid, time step, initialization, and stopping criterion.  Let \(\bar i\)
denote the direction opposite to an incoming
wall direction \(i\), let
\(\mat R_b=\diag(1,-1,-1)\), and write the prescribed wall state as
\(w_b=(\rho_b,h\rho_bu_{b,1},h\rho_bu_{b,2})^T\).  At the cell center \(x_f\)
adjacent to the boundary, we use the midpoint wall closure for vector-valued
populations introduced by Zhao, Zhang, and Yong~\cite{ZhaoZhangYong2020Boundary}:
\begin{equation}
 f_i^{n+1}(x_f)=\mat R_b f_{\bar i}^{*,n}(x_f)
 +M_i(w_b)-\mat R_b M_{\bar i}(w_b).
 \label{eq:channel_midpoint_boundary}
\end{equation}
Here the superscript \(*\) denotes the post-collision, pre-streaming state,
and all channel calculations set \(\rho_b=1\).
The cell centers are \(y_j=(j+1/2)h\), \(j=0,\ldots,n_g-1\), with
\(h=1/n_g\), so the walls lie midway between the
first or last cell center and its missing neighbor.  The boundary formula is
inherited from Ref.~\cite{ZhaoZhangYong2020Boundary}; the present tests examine
its coupling to the MRT collision and the resulting rate dependence.  For a
stationary wall, the equilibrium correction in
Eq.~\eqref{eq:channel_midpoint_boundary} vanishes.
After moment-space collision, the uniform streamwise body force is applied
only to the rest population as
\begin{equation}
 f_5^*\leftarrow f_5^*+\alpha h^3(0,G,0)^T.
 \label{eq:channel_body_force}
\end{equation}
On cell centers in $0<y<1$, the exact family is
\begin{equation}
u_x(y)=\frac{G}{2\nu}y(1-y)+U_{\rm top}y,\qquad u_y=0,
\quad Q=\frac{G}{12\nu}+\frac{U_{\rm top}}2.
\end{equation}
These solutions have \(\rho=1\) and \(P=0\).
We consider Poiseuille flow $(G,U_{\rm top})=(0.8\nu,0)$, Couette flow
$(G,U_{\rm top})=(0,0.05)$, and Couette--Poiseuille flow
$(G,U_{\rm top})=(0.8\nu,0.05)$.  Only steady, spatially uniform forcing is
considered.  In this unidirectional setting, the Couette--Poiseuille
calculation primarily checks superposition; it does not test additional
nonlinear effects.

Each calculation is initialized with the corresponding continuum velocity
profile.  The numerical solution is sampled at intervals $\Delta T=0.5$ and
stopped at the first pair of successive sampling times satisfying
\begin{equation}
 \frac{\|\vect u_h(T+\Delta T)-\vect u_h(T)\|_{2,h}}
 {\max\{\|\vect u_h(T+\Delta T)\|_{2,h},10^{-14}\}}<10^{-10}.
 \label{eq:channel_stopping_criterion}
\end{equation}
All runs reported in Table~\ref{tab:boundary} satisfy this criterion at or
before the maximum time $T=70$.
The initialized pure Couette profile is preserved to roundoff and is therefore
omitted from the error tables.
Let \(\bar u_{x,h}\) denote the streamwise average of the numerical velocity
and set $Q_h=\langle\bar u_{x,h}\rangle_h$.  The boundary errors are
\(E_u=\|\bar u_{x,h}-u_x\|_{2,h}/\|u_x\|_{2,h}\) and
\(E_Q=|Q_h-Q|/|Q|\).  Table~\ref{tab:boundary} reports these errors.

\begin{table*}[t]
\centering\footnotesize
\caption{Channel errors on the finest grid and orders fitted over
\(n_g=(8,12,16,24)\).  F-MRT and F-TRT use the Fourier-selected D2 channel
rates in Table~\ref{tab:selected-rates}; no channel calibration is used here.
Values at every grid level are given in
Section~S7 and Table~S6 of the Supplementary Material.}
\label{tab:boundary}
\begin{tabular}{llrrrr}
\toprule
Flow & rate set & \(E_u^{\rm fine}\) & \(r_u\) &
\(E_Q^{\rm fine}\) & \(r_Q\)\\
\midrule
Poiseuille & F-MRT
 & \(8.043\times10^{-4}\) & 1.9999 & \(1.303\times10^{-5}\) & 2.0000\\
& F-TRT
 & \(6.079\times10^{-4}\) & 1.9999 & \(1.534\times10^{-3}\) & 2.0000\\
& SRT
 & \(2.374\times10^{-3}\) & 1.9999 & \(1.733\times10^{-3}\) & 2.0000\\
\addlinespace
Couette--Poiseuille & F-MRT
 & \(6.026\times10^{-4}\) & 1.9990 & \(9.480\times10^{-6}\) & 1.9999\\
& F-TRT
 & \(4.554\times10^{-4}\) & 1.9990 & \(1.116\times10^{-3}\) & 2.0000\\
& SRT
 & \(1.779\times10^{-3}\) & 1.9990 & \(1.260\times10^{-3}\) & 2.0000\\
\bottomrule
\end{tabular}
\end{table*}

\begin{figure*}[t]
\centering
\includegraphics[width=0.96\textwidth]{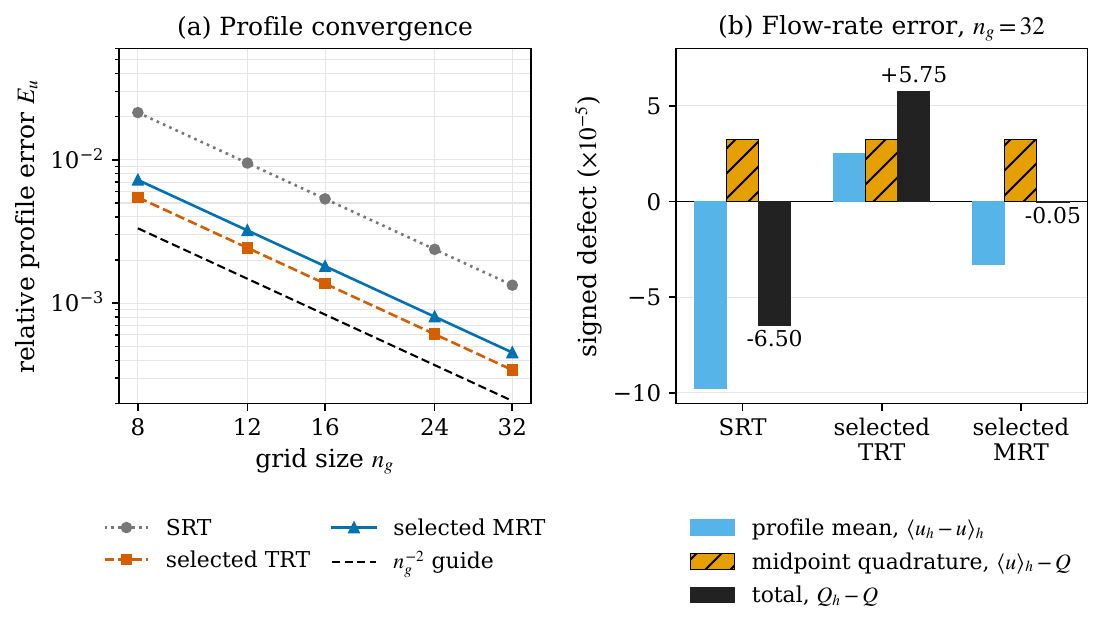}
\caption{Poiseuille flow with the midpoint wall closure of
Ref.~\cite{ZhaoZhangYong2020Boundary}.  The ``selected'' MRT and TRT
legend entries denote the Fourier-selected channel rates.
Panel (a) shows convergence of the velocity profile error over five grids; the
black dashed segment indicates $O(h^2)$.  Panel (b) shows the discrete mean
of the profile error, the common midpoint quadrature error, and their signed
sum in Eq.~\eqref{eq:flow_rate_defect_split} at $n_g=32$.  Their near
cancellation for Fourier-selected MRT explains its small flow-rate error, despite
its larger profile error than Fourier-selected TRT.}
\label{fig:midpoint-boundary}
\end{figure*}

Table~\ref{tab:boundary} gives fitted orders close to two for both the
profile and the flow rate; Figure~\ref{fig:midpoint-boundary}(a) shows the
Poiseuille profile-error curves, including an additional $n_g=32$ calculation.
On the finest grid in Table~\ref{tab:boundary}, the Fourier-selected MRT
profile error is about 32\% larger than that of Fourier-selected TRT, whereas its
error in flow rate is about 99.2\% smaller.  To explain this difference, we
decompose the flow-rate error as
\begin{equation}
 Q_h-Q=\langle\bar u_{x,h}-u_x\rangle_h
       +\bigl(\langle u_x\rangle_h-Q\bigr),
 \label{eq:flow_rate_defect_split}
\end{equation}
on an additional $n_g=32$ grid.  The discrete mean of the profile error and the
common midpoint quadrature error nearly cancel for Fourier-selected MRT:
their sum is only 0.85\% of the defect in the Fourier-selected TRT flow rate
despite the larger Fourier-selected MRT profile
error.  The smaller flow-rate error therefore does not imply a uniformly
more accurate velocity field.  Figure~\ref{fig:midpoint-boundary}(b)
shows the two signed components and their sum.  Section~S7 and Table~S6 of
the Supplementary Material give the profile and flow-rate errors on all
four grids.  The supplementary code computes the nodal profiles,
the decomposition at $n_g=32$, and the conservation diagnostics.

\paragraph{Channel calibration and fixed TRT references}
To assess the dependence of the channel results on the TRT rate, we also
considered the choice with complementary rates $s^++s_\nu=2$, the classical
reference for scalar LBM
$\Lambda=(s_\nu^{-1}-1/2)((s^+)^{-1}-1/2)=3/16$~\cite{Ginzburg2005}, and two
rates obtained by channel calibration, with
$s^+=s_\sigma=s_\varepsilon$.  The value $\Lambda=3/16$ is included only as
a reference because it is not known to be exact for the present vector-valued
scheme, forcing discretization, and wall update.
This calibration is independent of the TGV calibration above.  The two
channel-calibrated TRT rates were fitted only to the $n_g=8$ Poiseuille case
by minimizing
$\max(E_u,E_Q)$ over
$0.05\le s^+\le1.95$.  A 21-point initial grid followed by seven
refinements evaluated 46 distinct rates.  The $n_g=12,16,24$ cases and all
Couette--Poiseuille cases were excluded from this fit.  Using the boundary
errors alone gives
$s^+=0.114570$, but the collision--streaming Fourier symbol then has a
locally refined spectral radius of $1.008271>1$.  Requiring the prescribed
spectral condition, based on finite sampling and local optimization, instead
gives $s^+=0.083398$.  Both are channel-calibrated choices: the objective
remains $\max(E_u,E_Q)$, and the spectral condition restricts the feasible
rates without introducing the Fourier accuracy objective $E_T$.
These additional computations use the stopping
tolerance in Eq.~\eqref{eq:channel_stopping_criterion} and a maximum time
$T=120$; all
reach the tolerance by $T=72.5$.

\begin{table*}[t]
\centering\scriptsize
\caption{Channel errors at $n_g=24$ for cases excluded from calibration, using
F-MRT, F-TRT, and additional TRT references and channel-calibrated rates.
The latter minimize $\max(E_u,E_Q)$ only at $n_g=8$ for Poiseuille flow,
with or without the spectral constraint.  They are independent of the
TGV-calibrated rates.  In each error column, the upper and lower entries
are $E_u$ and $E_Q$.  ``Spectral test'' indicates whether the rate passes
the numerical test in Section~\ref{sec:algorithm}.}
\label{tab:boundary-controls}
\setlength{\tabcolsep}{4pt}
\newcommand{\controlcell}[2]{\shortstack{$#1$\\[-1pt]$#2$}}
\begin{tabular}{@{}lcccc@{}}
\toprule
rate choice & $(s_\sigma,s_\varepsilon)$ & Spectral test & Poiseuille &
Couette--Poiseuille\\
\midrule
F-MRT & $(0.072884,1.950000)$ & pass
 & \controlcell{8.043\times10^{-4}}{1.303\times10^{-5}}
 & \controlcell{6.026\times10^{-4}}{9.480\times10^{-6}}\\
F-TRT & $(0.076586,0.076586)$ & pass
 & \controlcell{6.079\times10^{-4}}{1.534\times10^{-3}}
 & \controlcell{4.554\times10^{-4}}{1.116\times10^{-3}}\\
\shortstack[l]{Channel-calibrated TRT\\with spectral constraint} & $(0.083398,0.083398)$ & pass
 & \controlcell{3.544\times10^{-4}}{1.256\times10^{-3}}
 & \controlcell{2.655\times10^{-4}}{9.136\times10^{-4}}\\
\shortstack[l]{Channel-calibrated TRT\\without spectral constraint} & $(0.114570,0.114570)$ & fail
 & \controlcell{4.212\times10^{-4}}{4.067\times10^{-4}}
 & \controlcell{3.155\times10^{-4}}{2.958\times10^{-4}}\\
C-TRT & $(0.048780,0.048780)$ & pass
 & \controlcell{2.377\times10^{-3}}{3.472\times10^{-3}}
 & \controlcell{1.781\times10^{-3}}{2.525\times10^{-3}}\\
$\Lambda=3/16$ TRT & $(0.064516,0.064516)$ & pass
 & \controlcell{1.189\times10^{-3}}{2.170\times10^{-3}}
 & \controlcell{8.905\times10^{-4}}{1.578\times10^{-3}}\\
SRT & $(1.951220,1.951220)$ & pass
 & \controlcell{2.374\times10^{-3}}{1.733\times10^{-3}}
 & \controlcell{1.779\times10^{-3}}{1.260\times10^{-3}}\\
\bottomrule
\end{tabular}
\end{table*}

Table~\ref{tab:boundary-controls} compares the two flows at $n_g=24$, neither
of which was used for channel calibration.  Channel-calibrated TRT without the spectral
constraint gives a smaller profile error and a smaller $\max(E_u,E_Q)$ than
Fourier-selected MRT, but fails the spectral test.  Among the rates passing this
test, channel-calibrated TRT gives the smaller profile error, whereas Fourier-selected MRT gives
the smaller flow-rate error and the smaller $\max(E_u,E_Q)$.
No TRT rate satisfying the spectral constraint in the finite search attained
the Fourier-selected MRT accuracy in flow rate.  The supplementary code
includes the channel calibration and grid-refinement calculations used in
this comparison.

\section{Conclusions}
\label{sec:conclusions}

We have constructed D2N5 and D3N7 vectorial MRT schemes with a common
velocity-moment operator for all physical components and independent flux,
stress, and energy relaxation rates.  The stated rate
equalities recover their TRT and SRT subfamilies.  Maxwell iteration gives
formal second-order consistency for the unforced interior schemes under
diffusive scaling.  Only the flux rate determines the viscosity; for D2N5,
the explicit conserved expansion identifies the first contribution of the
stress and energy rates to the $O(h^2)$ macroscopic truncation error.
A lower triangular change of variables yields a symmetrizer for the full
collision Jacobian at rest.  Weighted moment orthogonality also gives an exact
identity for the nonlinear kinetic residual under local collision.

Fourier selection uses a common accuracy and damping criterion for MRT and
TRT.  In all three parameter cases, the selected MRT rates give a smaller
transverse eigenvalue error and a smaller locally maximized high-frequency
spectral radius than the selected TRT rates.  The resulting spectral
comparison differs from the finite-time velocity-error comparison, which
also depends on modal projection and kinetic transients.

The exact Taylor--Green and ABC solutions and the oblique transverse plane
waves exhibit second-order convergence.  For the nonlinear D3 Taylor--Green
test, coarse-grid calibration reduces the MRT velocity and vorticity errors
by 37.4\% and 33.7\%, respectively, relative to similarly calibrated TRT on
the finest grid.  The two fixed rate pairs retain approximately second-order
convergence on the finer grids, and the MRT errors remain smaller at the
higher amplitude.  The tested implementations have comparable costs per
lattice site and time step; MRT therefore attains the prescribed accuracy
with less estimated simulation time on these grids.  This timing comparison
excludes reference calculations and rate calibration, for which the MRT
search is more expensive.

Complementary TRT
has the smallest maximum velocity error over the four additional plane waves
at $n_g=32$, and SRT has the smallest field errors in the strong double shear
layer test.  The channel calculations exhibit approximately second-order
convergence with the previously developed midpoint wall scheme.  There,
Fourier-selected MRT gives a smaller flow-rate error but a larger profile error than
Fourier-selected TRT; the small flow-rate error partly results from cancellation of
two $O(h^2)$ contributions.

The parameter comparisons are conditional on the finite search sets; in
particular, the Fourier-selected channel energy rate lies at the upper
bound of the prescribed interval.  Stability and convergence of the fully
discrete initial-boundary-value problem require further analysis of the
collision--streaming operator and its coupling with boundary conditions.

\section*{Acknowledgments}
The authors thank Professor Weifeng Zhao for helpful discussions.  This work
was supported by the Beijing Natural Science Foundation (Grant
No.\ JR25003) and the Beijing Outstanding Young Scientist Program
(No.\ JWZQ20240101027).

\bigskip
\noindent\begin{minipage}{\textwidth}
\section*{CRediT authorship contribution statement}
\textbf{Yuyang Xu:} Software, Validation, Writing -- review \& editing.\\
\textbf{Yuxuan Zhang:} Writing -- review \& editing.\\
\textbf{Jin Zhao:} Conceptualization, Methodology, Software, Formal analysis,
Writing -- original draft, Writing -- review \& editing.
\end{minipage}

\section*{Declaration of competing interest}
The authors declare that there are no known competing financial interests or
personal relationships that could have appeared to influence the work
reported in this paper.

\section*{Data availability}
Numerical results are reported in this article and the Supplementary Material.
Standalone Python code, fixed input parameters, and implementation tests are
provided in the supplementary code package.  Reference solutions are computed
when the code is run.

\bigskip
\noindent\begin{minipage}{\textwidth}
\section*{Declaration of generative AI and AI-assisted technologies in the manuscript preparation process}
The authors used OpenAI Codex to assist with manuscript preparation.
The authors reviewed and revised the resulting material and take
responsibility for the content of this article.
\end{minipage}

\begingroup
\small
\hbadness=2000
\bibliographystyle{elsarticle-num}
\bibliography{AML-jcp}

@article{Chen1998,
  author  = {Chen, Shiyi and Doolen, Gary D.},
  title   = {Lattice {B}oltzmann method for fluid flows},
  journal = {Annual Review of Fluid Mechanics},
  volume  = {30},
  pages   = {329--364},
  year    = {1998},
  doi     = {10.1146/annurev.fluid.30.1.329}
}

@book{Succi2001,
  author    = {Succi, Sauro},
  title     = {The Lattice {B}oltzmann Equation for Fluid Dynamics and Beyond},
  publisher = {Oxford University Press},
  year      = {2001},
  address   = {Oxford},
  doi       = {10.1093/oso/9780198503989.001.0001}
}

@book{Guo2013,
  author    = {Guo, Zhaoli and Shu, Chang},
  title     = {Lattice {B}oltzmann Method and Its Applications in Engineering},
  publisher = {World Scientific},
  year      = {2013},
  address   = {Singapore},
  doi       = {10.1142/8806}
}

@article{Bhatnagar1954,
  author  = {Bhatnagar, P. L. and Gross, E. P. and Krook, M.},
  title   = {A model for collision processes in gases. {I}. {S}mall amplitude processes in charged and neutral one-component systems},
  journal = {Physical Review},
  volume  = {94},
  pages   = {511--525},
  year    = {1954},
  doi     = {10.1103/PhysRev.94.511}
}

@incollection{DHumieres1992,
  author    = {d'Humi\`{e}res, Dominique},
  title     = {Generalized lattice-{B}oltzmann equations},
  booktitle = {Rarefied Gas Dynamics: Theory and Simulations},
  editor    = {Shizgal, Bernie D. and Weaver, Donald P.},
  series    = {Progress in Astronautics and Aeronautics},
  volume    = {159},
  pages     = {450--458},
  publisher = {American Institute of Aeronautics and Astronautics},
  year      = {1994},
  doi       = {10.2514/5.9781600866319.0450.0458}
}

@article{Lallemand2000,
  author  = {Lallemand, Pierre and Luo, Li-Shi},
  title   = {Theory of the lattice {B}oltzmann method: Dispersion, dissipation, isotropy, {G}alilean invariance, and stability},
  journal = {Physical Review E},
  volume  = {61},
  pages   = {6546--6562},
  year    = {2000},
  doi     = {10.1103/PhysRevE.61.6546}
}

@article{Ginzburg2008,
  author  = {Ginzburg, Irina and Verhaeghe, Frederik and d'Humi\`{e}res, Dominique},
  title   = {Two-relaxation-time lattice {B}oltzmann scheme: {A}bout parametrization, velocity, pressure and mixed boundary conditions},
  journal = {Communications in Computational Physics},
  volume  = {3},
  number  = {2},
  pages   = {427--478},
  year    = {2008},
  doi     = {10.4208/cicp.2008.v3.p427}
}

@article{Ginzburg2005,
  author  = {Ginzburg, Irina and d'Humi\`{e}res, Dominique},
  title   = {Multireflection boundary conditions for lattice {B}oltzmann models},
  journal = {Physical Review E},
  volume  = {68},
  pages   = {066614},
  year    = {2003},
  doi     = {10.1103/PhysRevE.68.066614}
}

@article{Dellar2002,
  author  = {Dellar, Paul J.},
  title   = {Nonhydrodynamic modes and a priori construction of shallow water lattice {B}oltzmann equations},
  journal = {Physical Review E},
  volume  = {65},
  number  = {3},
  pages   = {036309},
  year    = {2002},
  doi     = {10.1103/PhysRevE.65.036309}
}

@article{Luo2011,
  author  = {Luo, Li-Shi and Liao, Wei and Chen, Xiaopeng and Peng, Yan and Zhang, Wei},
  title   = {Numerics of the lattice {B}oltzmann method: {E}ffects of collision models on the lattice {B}oltzmann simulations},
  journal = {Physical Review E},
  volume  = {83},
  pages   = {056710},
  year    = {2011},
  doi     = {10.1103/PhysRevE.83.056710}
}

@book{Kruger2017,
  author    = {Kr\"{u}ger, Timm and Kusumaatmaja, Halim and Kuzmin, Alexandr and Shardt, Orest and Silva, Gon{\c c}alo and Viggen, Erlend Magnus},
  title     = {The Lattice {B}oltzmann Method: Principles and Practice},
  series    = {Graduate Texts in Physics},
  publisher = {Springer},
  address   = {Cham},
  year      = {2017},
  isbn      = {978-3-319-44647-9},
  doi       = {10.1007/978-3-319-44649-3}
}

@article{Geier2006,
  author  = {Geier, Martin and Greiner, Andreas and Korvink, Jan G.},
  title   = {Cascaded digital lattice {B}oltzmann automata for high {R}eynolds number flow},
  journal = {Physical Review E},
  volume  = {73},
  number  = {6},
  pages   = {066705},
  year    = {2006},
  doi     = {10.1103/PhysRevE.73.066705}
}

@article{Geier2015,
  author  = {Geier, Martin and Sch\"{o}nherr, Martin and Pasquali, Andrea and Krafczyk, Manfred},
  title   = {The cumulant lattice {B}oltzmann equation in three dimensions: Theory and validation},
  journal = {Computers \& Mathematics with Applications},
  volume  = {70},
  number  = {4},
  pages   = {507--547},
  year    = {2015},
  doi     = {10.1016/j.camwa.2015.05.001}
}

@article{Karlin1999,
  author  = {Karlin, Iliya V. and Ferrante, Alessandro and \"{O}ttinger, Hans Christian},
  title   = {Perfect entropy functions of the lattice {B}oltzmann method},
  journal = {Europhysics Letters},
  volume  = {47},
  number  = {2},
  pages   = {182--188},
  year    = {1999},
  doi     = {10.1209/epl/i1999-00370-1}
}

@article{Karlin2014,
  author  = {Karlin, Iliya V. and B\"{o}sch, Fabian and Chikatamarla, Shyam S.},
  title   = {Gibbs' principle for the lattice-kinetic theory of fluid dynamics},
  journal = {Physical Review E},
  volume  = {90},
  number  = {3},
  pages   = {031302},
  year    = {2014},
  doi     = {10.1103/PhysRevE.90.031302}
}

@article{Suga2015,
  author  = {Suga, Kazuhiko and Kuwata, Yusuke and Takashima, Kazuki and Chikasue, Rina},
  title   = {A {D3Q27} multiple-relaxation-time lattice {B}oltzmann method for turbulent flows},
  journal = {Computers \& Mathematics with Applications},
  volume  = {69},
  number  = {6},
  pages   = {518--529},
  year    = {2015},
  doi     = {10.1016/j.camwa.2015.01.010}
}

@article{Coreixas2019,
  author  = {Coreixas, Christophe and Chopard, Bastien and Latt, Jonas},
  title   = {Comprehensive comparison of collision models in the lattice {B}oltzmann framework: Theoretical investigations},
  journal = {Physical Review E},
  volume  = {100},
  pages   = {033305},
  year    = {2019},
  doi     = {10.1103/PhysRevE.100.033305}
}

@article{Wissocq2019,
  author  = {Wissocq, Gauthier and Sagaut, Pierre and Boussuge, Jean-Fran\c{c}ois},
  title   = {An extended spectral analysis of the lattice {B}oltzmann method: Modal interactions and stability issues},
  journal = {Journal of Computational Physics},
  volume  = {380},
  pages   = {311--333},
  year    = {2019},
  doi     = {10.1016/j.jcp.2018.12.015}
}

@article{Hosseini2020,
  author  = {Hosseini, Seyed Ali and Coreixas, Christophe and Darabiha, Nasser and Th\'{e}venin, Dominique},
  title   = {Stability of the lattice kinetic scheme and choice of the free relaxation parameter},
  journal = {Physical Review E},
  volume  = {99},
  number  = {6},
  pages   = {063305},
  year    = {2019},
  doi     = {10.1103/PhysRevE.99.063305}
}

@article{Taylor1937,
  author  = {Taylor, Geoffrey Ingram and Green, Albert Edward},
  title   = {Mechanism of the production of small eddies from large ones},
  journal = {Proceedings of the Royal Society of London. Series A - Mathematical and Physical Sciences},
  volume  = {158},
  number  = {895},
  pages   = {499--521},
  year    = {1937},
  doi     = {10.1098/rspa.1937.0036}
}

@article{Dellar2025VectorProjection,
  author  = {Dellar, Paul J.},
  title   = {A projection method to symmetrize the momentum flux in a vector lattice {B}oltzmann formulation of hydrodynamics},
  journal = {Physics of Fluids},
  volume  = {37},
  number  = {6},
  pages   = {063606},
  year    = {2025},
  doi     = {10.1063/5.0268328}
}

@misc{Zhao2026TRTPower,
  author        = {Zhao, Jin},
  title         = {Mesh-uniform power stability of two-relaxation-time vector lattice {B}oltzmann schemes with reversible boundaries},
  year          = {2026},
  eprint        = {2608.30253},
  archiveprefix = {arXiv},
  primaryclass  = {math.NA},
  doi           = {10.48550/arXiv.2608.30253}
}

@article{zhao2017maxwell,
  author  = {Zhao, Weifeng and Yong, Wen-An},
  title   = {Maxwell iteration for the lattice {B}oltzmann method with diffusive scaling},
  journal = {Physical Review E},
  volume  = {95},
  number  = {3},
  pages   = {033311},
  year    = {2017},
  doi     = {10.1103/PhysRevE.95.033311}
}

@article{CarforaNatalini2008,
  author  = {Carfora, Maria Francesca and Natalini, Roberto},
  title   = {A discrete kinetic approximation for the incompressible {N}avier--{S}tokes equations},
  journal = {ESAIM: Mathematical Modelling and Numerical Analysis},
  volume  = {42},
  number  = {1},
  pages   = {93--112},
  year    = {2008},
  doi     = {10.1051/m2an:2007055}
}

@article{BouchutEtAl2018BGKFVS,
  author  = {Bouchut, Fran\c{c}ois and Jobic, Yann and Natalini, Roberto and Occelli, Ren\'{e} and Pavan, Vincent},
  title   = {Second-order entropy satisfying {BGK-FVS} schemes for incompressible {N}avier--{S}tokes equations},
  journal = {The SMAI Journal of Computational Mathematics},
  volume  = {4},
  pages   = {1--56},
  year    = {2018},
  doi     = {10.5802/smai-jcm.28}
}

@article{ZhaoZhangYong2020Approximation,
  author  = {Zhao, Jin and Zhang, Zhimin and Yong, Wen-An},
  title   = {Approximation of the multi-dimensional incompressible {N}avier--{S}tokes equations by discrete-velocity vector-{BGK} models},
  journal = {Journal of Mathematical Analysis and Applications},
  volume  = {486},
  number  = {2},
  pages   = {123901},
  year    = {2020},
  doi     = {10.1016/j.jmaa.2020.123901}
}

@article{ZhaoZhangYong2020Boundary,
  author  = {Zhao, Jin and Zhang, Zhimin and Yong, Wen-An},
  title   = {Vector-Type Boundary Schemes for the Lattice {B}oltzmann Method Based on Vector-{BGK} Models},
  journal = {SIAM Journal on Scientific Computing},
  volume  = {42},
  number  = {5},
  pages   = {B1250--B1270},
  year    = {2020},
  doi     = {10.1137/19M1308542}
}

@article{Zhao2021VectorBGK,
  author  = {Zhao, Jin},
  title   = {Discrete-velocity vector-{BGK} models based numerical methods for the incompressible {N}avier--{S}tokes equations},
  journal = {Communications in Computational Physics},
  volume  = {29},
  number  = {2},
  pages   = {420--444},
  year    = {2021},
  doi     = {10.4208/cicp.OA-2019-0192}
}

@article{ZhaoYong2020Weighted,
  author  = {Zhao, Weifeng and Yong, Wen-An},
  title   = {Weighted {L}$^2$-stability of a discrete kinetic approximation for the incompressible {N}avier--{S}tokes equations on bounded domains},
  journal = {Journal of Computational and Applied Mathematics},
  volume  = {376},
  pages   = {112820},
  year    = {2020},
  doi     = {10.1016/j.cam.2020.112820}
}

@article{ZhaoYong2024VFDLBM,
  author  = {Zhao, Jin and Yong, Wen-An},
  title   = {Vectorial finite-difference-based lattice {B}oltzmann method: Consistency, boundary schemes and stability analysis},
  journal = {Journal of Computational and Applied Mathematics},
  volume  = {441},
  pages   = {115677},
  year    = {2024},
  doi     = {10.1016/j.cam.2023.115677}
}

@article{YangZhaoLin2024Automatic,
  author  = {Yang, Jianbin and Zhao, Weifeng and Lin, Ping},
  title   = {An automatic approach for the stability analysis of multi-relaxation-time lattice {B}oltzmann models},
  journal = {Journal of Computational Physics},
  volume  = {519},
  pages   = {113432},
  year    = {2024},
  doi     = {10.1016/j.jcp.2024.113432}
}

@article{Dellar2003IncompressibleLimits,
  author  = {Dellar, Paul J.},
  title   = {Incompressible limits of lattice {B}oltzmann equations using multiple relaxation times},
  journal = {Journal of Computational Physics},
  volume  = {190},
  number  = {2},
  pages   = {351--370},
  year    = {2003},
  doi     = {10.1016/S0021-9991(03)00279-1}
}

@article{ZhaoZhang2020D3N7,
  author  = {Zhao, Jin and Zhang, Zhimin},
  title   = {Seven-velocity three-dimensional vectorial lattice {B}oltzmann method including various types of approximations to the pressure and two-parameterized second-order boundary treatments},
  journal = {Computers \& Mathematics with Applications},
  volume  = {80},
  number  = {12},
  pages   = {2764--2779},
  year    = {2020},
  doi     = {10.1016/j.camwa.2020.10.001}
}

@misc{AregbaDriolletEtAl2026Vectorial,
  author        = {Aregba-Driollet, Denise and Bellotti, Thomas and Natalini, Roberto and Tenna, Tommaso},
  title         = {A vectorial lattice {B}oltzmann scheme for the incompressible {N}avier--{S}tokes equations},
  year          = {2026},
  eprint        = {2607.24059},
  archiveprefix = {arXiv},
  primaryclass  = {math.NA},
  doi           = {10.48550/arXiv.2607.24059}
}

@article{Dubois2014Vectorial,
  author  = {Dubois, Fran\c{c}ois},
  title   = {Simulation of strong nonlinear waves with vectorial lattice {B}oltzmann schemes},
  journal = {International Journal of Modern Physics C},
  volume  = {25},
  number  = {12},
  pages   = {1441014},
  year    = {2014},
  doi     = {10.1142/S0129183114410149}
}

@article{XuSagaut2011Optimal,
  author  = {Xu, Hui and Sagaut, Pierre},
  title   = {Optimal low-dispersion low-dissipation {LBM} schemes for computational aeroacoustics},
  journal = {Journal of Computational Physics},
  volume  = {230},
  number  = {13},
  pages   = {5353--5382},
  year    = {2011},
  doi     = {10.1016/j.jcp.2011.03.040}
}

@article{XuMalaspinasSagaut2012Sensitivity,
  author  = {Xu, Hui and Malaspinas, Orestis and Sagaut, Pierre},
  title   = {Sensitivity analysis and determination of free relaxation parameters for the weakly-compressible {MRT--LBM} schemes},
  journal = {Journal of Computational Physics},
  volume  = {231},
  number  = {21},
  pages   = {7335--7367},
  year    = {2012},
  doi     = {10.1016/j.jcp.2012.07.005}
}

@article{BellColellaGlaz1989,
  author  = {Bell, John B. and Colella, Phillip and Glaz, Harland M.},
  title   = {A second-order projection method for the incompressible {N}avier--{S}tokes equations},
  journal = {Journal of Computational Physics},
  volume  = {85},
  number  = {2},
  pages   = {257--283},
  year    = {1989},
  doi     = {10.1016/0021-9991(89)90151-4}
}
\endgroup

\end{document}